\documentclass[12pt]{amsart}

\usepackage[utf8]{inputenc}
\usepackage[english]{babel}
\usepackage{amsmath}
\usepackage{amsfonts}
\usepackage{amssymb}
\usepackage{graphicx}
\usepackage{mathtools}
\usepackage{amsthm}
\newtheorem*{Thm}{Theorem}

\newtheorem{thm}{Theorem}[section]
\newtheorem{prop}[thm]{Proposition}
\newtheorem{lem}[thm]{Lemma}

\newtheorem{rem}[thm]{Remark}

\theoremstyle{definition}
\newtheorem{defi}{Definition}[section]
\newcommand{\R}{\mathbb{R}}
\renewcommand{\P}{\mathbb{P}}

\newcommand{\Z}{\mathbb{Z}}
\newcommand{\C}{\mathbb{C}}

\newcommand{\Conf}{\mathrm{Conf}}

\usepackage{fullpage}

\usepackage{hyperref}

\title{Explicit description of Christoffel deformations and Palm measures of the Plancherel measure, the $z$-measures and the Gamma process}

\author[P. Lazag]{Pierre Lazag}

\begin{document}
\begin{abstract}Christoffel deformation of a measure on the real line consists of multipying this measure by a squared polynomial having its roots in $\R$. We introduce Christoffel deformations of discrete orthogonal polynomial ensembles by considering the Christoffel deformations of the underlying measure, and prove that this construction extends to more general point processes describing distribution on partitions: the poissonized Plancherel measure and the $z$-measures. These deformations contain the theory of Palm measures, and for example, explicit formulas for the Palm measures of the poissonized Plancherel measure provide a description of the TASEP with initial wedge condition with frozen particles. We also obtain new formulas for Palm measures of the $z$-measures. The extension to the Plancherel measure is obtained via a limit transition from the Charlier ensemble, while the extension to the $z$-measures follows from an analytic continuation argument. A limit procedure starting from the non-degenerate $z$-measures leads to a deformation of the Gamma process introduced by Borodin and Olshanski.
\end{abstract}
\maketitle
\scriptsize{\address{CNRS, UMR 7373, \'Ecole centrale de Marseille, Institut de Math\'ematiques de Marseille, Aix-Marseille Universit\'e, Marseille, France}

\address{SISSA, via Bonomea 265, 34136, Trieste, Italy}

\address{INFN, Section of Trieste}

\email pierrelazag@hotmail.fr }
\\

\normalsize
Keywords: determinantal point processes, Palm measures, Christoffel deformations, Plancherel measures, $z$-measures, Gamma process

\section{Introduction}
\subsection{Random partitions as point processes}
The aim of this article is to study deformations of discrete determinantal point processes describing distributions or limit distributions on partitions, namely the poissonized Plancherel measure, the $z$-measures and its limiting Gamma process (see e.g. \cite{boo}, \cite{johanssonplancherel}, \cite{infinitewedge}, \cite{bodist}, \cite{bogamma} and \cite{bomarkov}). Recall that a partition $\lambda=(\lambda_1 \geq \lambda_2 \geq ...)$ is a finite non-increasing sequence of non-negative integers, and that it can be identified with a Young diagram (\cite{macdonald}). The set of all partitions will be denoted by $\mathbb{Y}$. As it is usual, for a partition $\lambda \in \mathbb{Y}$, we define its length $l(\lambda)$ as being the index of its last non-zero entry  
\begin{align*}l(\lambda):=\max \{i, \hspace{0.1cm} \lambda_i \neq 0 , \hspace{0.1cm} i=1,2,... \},
\end{align*}
and denote its size by $|\lambda|$ 
\begin{align*}
|\lambda|:= \sum_{i \geq 1} \lambda_i.
\end{align*}
A Young diagram $\lambda \in \mathbb{Y}$ can be seen as a subset of $\mathbb{Z}$ via the map $\lambda \mapsto \{ \lambda_i - i , i=1,2,... \}$ or some shifted version of it. More precisely, when we consider partitions $\lambda \in \mathbb{Y}$ with fixed length $N$, it is convenient to consider the finite set of positive integers 
\begin{align*}
\{\lambda_i -i+N, \hspace{0.1cm} i=1,...,N\},
\end{align*}
while when the lengths may vary, it is convenient to consider the infinite sets 
\begin{align*}
\{ \lambda_i - i , \hspace{0.1cm} i=1,2,... \} \text{ or }\{ \lambda_i - i +1/2, \hspace{0.1cm} i=1,2,... \},
\end{align*}
the latter providing a more symmetric picture. In any case, the pushforward of a probability measure on $\mathbb{Y}$ by one of these maps is thus a point process on a discrete subset of $\mathbb{R}$. In the cases we deal with, namely the poissonized Plancherel measure, the $z$-measures and its limiting process with the Gamma kernel, these point processes are determinantal point processes, see definition \ref{defdpp} below.
\subsection{Determinantal point processes and Christoffel deformations}
A specific and important class of determinantal point processes are the so called \textit{orthogonal polynomial ensembles}, see e.g. the survey \cite{konig} and references therein. They are $N$-points random configurations for some fixed deterministic $N$ governed by the kernel of the orthogonal projection onto the first $N$ polynomials in the Hilbert space $L^2(\mathbb{R},\mu)$ for some measure $\mu$ with finite moments of all orders. This kernel is namely the $N$-th Christoffel-Darboux kernel of the polynomials orthogonal with respect to $\mu$. For example, the classical G.U.E. from random matrix theory is the orthogonal polynomial ensemble with the Gaussian weight on $\mathbb{R}$, and is thus described by the Hermite polynomials, see e.g. \cite{konig} and references therein.\\
\par
Multiplying the weight $\mu$ by a squared polynomial having its roots in $\R$ leads to the so-called Christoffel deformation of the system and an explicit description of the new family of orthogonal polynomials as well as the corresponding deformed point process.
\subsection{Background and motivations: characteristic polynomials and Palm measures}
We prove in this note that the multiplicative structure of Christoffel deformations of orthogonal polynomial ensembles carries on to more general determinantal point processes. Our method is essentially by approximation of the models in consideration by orthogonal polynomial ensembles. However, since from \cite{bufetovconditional} we know that any determinantal point process on the real line is approximated by orthogonal polynomial ensemble in a particular sense, we hope to establish the existence and the explicit descritption of Christoffel deformations of determinantal point processes in full generality in a future work; this method has already been fruitful in order to prove the Giambelli compatibility of determinantal point processes in \cite{bufetovlazag}.\\

Besides being a natural deformation of a point process, there are two directions in which Christoffel deformations have applications: characteristic polynomials and Palm measures. Let us also mention that Christoffel deformations have been considered in \cite{bufetovqiupickrell} in order to describe the ergodic decomposition of inifinite Hua-Pickrell measures.\\

The study of charcteristic polynomials of random matrices has connections with number theory, see for example \cite{keatingsnaith} or the more recent paper \cite{chhaibinajnudel}. In \cite{baikdeiftstrahov}, explicit formulas for averages of products and ratios of characteristic polynomials of random matrices are obtained by the formalism of Christoffel deformations.
\par
There are lots of studies of averages or more generally moments of characteristic polynomials in the continuous setting. For the discrete one, in \cite{giambelli}, the authors develop the notion of \textit{Giambelli compatible point processes} which provides a suitable formalism in the study of characteristic polynomials of general point processes. Among others, the authors show that the $z$-measures are Giambelli compatible, which means that characteristic polynomials of these processes are well defined and involve determinantal formulas.
\par
Rather than directly studying characteristic polynomials, from which determinantal formulas for the correlations can be obtained (see \cite{fyodorovstrahov1}, \cite{fyodorovstrahov2}, \cite{borodinstrahov}, and \cite{giambelli}), we go the reverse way and study Christoffel deformations of determinantal processes. Direct information on the characteristic polynomials of our processes of interest can be obtained from formulas from \cite{baikdeiftstrahov}, \cite{fyodorovstrahov1} or \cite{fyodorovstrahov2}. For instance, our results concerning the Charlier ensemble were used in \cite{bertolaruzza} in order to compute the characteristic polynomial of the discrete Bessel point process (which describes the poissonized Plancherel measure, see section \ref{sec3} below) and in connection with the Gromov-Witten theory of $\mathbb{P}^1$ (\cite{bertolaruzza}, Proposition 1.4). \\

On the other hand, the formalism of Christoffel deformations of orthogonal polynomial ensembles englobes the theory of Palm measures of the latter, see Proposition \ref{proppalm} below. Roughly speaking, the Palm measure of a given point process is the same point process conditioned to have particles at specific points. Our results thus provide new explicit and simple formulas for the Palm measures of the examples we are intersted in.
\par
An eloquent example is the one of the Totally Asymmetric Exclusion Process (TASEP) with initial wedge condition, which we now informally describe. We put independent exponential clock at each site $x \in \mathbb{Z}$, and at time $t=0$, we put particles on the whole semi-axis $\mathbb{Z}_{< 0}$. When a clock ring, the correponding particle moves to the site directly to the right of it, if this site is not occupied, and do not move else. A picture of this process at any time $t \geq 0$ is precisely the image of the poissonized Plancherel measure with parameter $t$ by the map 
\begin{align*}
\lambda \mapsto \{ \lambda_i - i , \hspace{0.1cm} i=1,2,... \},
\end{align*}
see for example \cite{borodingorin} and references therein. Theorem \ref{thmpalmbessel} thus gives a description of the TASEP with "frozen" sites, i.e. sites which must always be occupied, see Theorem \ref{thmpalmbessel} and remark \ref{remtasep} below.

In Theorem \ref{thmpalmzmeas}, we give new and compact formulas for the Palm measures of the $z$-measures.\\

Our results also connect with those of \cite{bufetovbessel} and \cite{bufetovolshanski} (see also \cite{bufetovquasisymmetries}). In these papers, the author shows that conditioning at different points leads to equivalent point processes and that the Radon-Nikodym derivatives are explicitly given by regularized multiplicative functionals, analogous to the one given by Christoffel deformations.\\
\par
The $z$-measures on partitions were introduced in \cite{bodist} in order to describe the harmonic analysis of the infinite symmetric group. One of their interset indeed lies in the fact that they constitute a coherent system of measures on the Young graph. However, we were not able to prove whether the Christoffel deformations of the $z$-measures we introduce constitute a coherent system or not, but we believe that the study of their characteristic polynomial and conditional measures is of interest.
\subsection{Main results}
As mentioned above, it is easy to define Christoffel deformations of orthogonal polynomial ensembles, but neither the Plancherel nor the $z$-measures are orthogonal polynomial ensembles. As for the characteristic polynomials, an obstacle is that, unlike for orthogonal polynomial ensembles, the number of points might be random or infinite. We still manage to give a meaning to Christoffel deformations of these processes, and also for the process with the Gamma kernel, and we now briefly explain our arguments and state our main results in an informal way.
\\
\par
On the one hand, we know from Johansson's paper \cite{johanssonplancherel} that the poissonized Plancherel measure is the limit of the Charlier ensemble, i.e. the orthogonal polynomial ensemble with the Poisson weight. The poissonized Plancherel measure is a determinantal point process described by Bessel functions, the discrete Bessel kernel (\cite{borodinbessel}). Refining asymptotic results of \cite{johanssonplancherel}, we prove the convergence of the Charlier polynomials towards the Bessel functions, see Lemma \ref{lemme1}, and obtain the following Theorem, see Theorems \ref{thm1}  and \ref{thmpalmbessel} for precise statements.
\begin{Thm} The Christoffel deformation of the Charlier ensemble converges to a determinantal point process described by discrete Wronskians of Bessel functions, the latter being a deformation of the poissonized Plancherel measure. The Palm measures of the poissonized Plancherel measure are described by similar kernels.
\end{Thm}
On the other hand, the $z$-measures are determinantal point processes described by specific Gauss hypergeometric functions. The analytic continuation argument presented in \cite{bomarkov} also allows us to describe the deformation of all $z$-measures, a particular case of them being the Meixner ensemble, and we obtain the following Theorem, see Theorems \ref{thmzmeas} and \ref{thmpalmzmeas} for precise statements.
\begin{Thm}The Christoffel deformations of the $z$-measures are determinantal point proccesses governed by a kernel involving discrete Wronskians of Gauss hypergeometric functions. Palm measures of the $z$-measures are governed by similar kernels.
\end{Thm}
Moreover, a somehow tricky but elementary asymptotic analysis, similar to the one performed in \cite{bogamma}, leads to a deformation of the process with the Gamma kernel, see Theorem \ref{thmgamma}.
\begin{Thm}The one-point deformations of the non-degenerate $z$-measures converge to a determinantal point process described by the Gamma function and its first derivative, which might be considered as a deformation of the process with the Gamma kernel.
\end{Thm}
We strongly believe that the approximation procedure we use in the last two Theorems could also be applied to the $zw$-measures, since the same argument of analytic continuation from an orthogonal polynomial ensemble holds, see \cite{bo-zw}.
\subsection{Organization of the paper}
The paper is organized as follows. In section \ref{sec2}, we collect definitions and basic properties of Christoffel deformations, orthogonal polnomial ensembles and determinantal point processes, and formally state that such deformations contain the description of their Palm measure. 
\par
In section \ref{sec3}, we introduce the Charlier ensemble and the Plancherel measure. We describe the Christoffel deformation of the Charlier ensemble and pass to the limit, in the same regime as in \cite{johanssonplancherel}. This leads to a deformation of the poissonized Plancherel measure, i.e. the process with the discrete Bessel kernel, see Theorem \ref{thm1}. We finally give the correlation kernel of the Palm measures of the discrete bessel point process.
\par
In section \ref{sec4}, we introduce the $z$-measures and define their Christoffel deformations. We describe the Christoffel deformations in the case where the parameters belong to the degenerate series (which is nothing but the Meixner ensemble) and prove that a similar determinantal formula also holds for the other series, see Theorem \ref{thmzmeas}. We further give a description of the Palm measures of the $z$-measures in Theorem \ref{thmpalmzmeas}. We then pass to the limit towards the deformed Gamma kernel in the case of the one point deformation, see Theorem \ref{thmgamma}.
\subsection*{Acknowledgments.}I am deeply grateful to Alexander Bufetov who posed the problem to me and for helpful discussions. I also would like to thank Alexander Boritchev and Pascal Hubert for helpful discussions and comments.\\

This project has received funding from the European Research Council (ERC) under the European Union's Horizon 2020 research and innovation programme (grant agreement No 647133 I-CHAOS).

The author acknowledges the support of the fellowship "Assegni di ricerca FSE SISSA
2019" from Fondo Sociale Europeo - Progetto SISSA OPERAZIONE 1 codice FP195673001.

The author acknowledges the financial support of the
project MMNLP (Mathematical Methods in Non Linear Physics) of the
INFN.

\subsection*{Conflict of interest}

The author declares that he has no conflict of interest.

\subsection*{Data availability}

Data sharing not applicable to this article as no datasets were generated or analysed during the current study.

\section{Christoffel deformations of orthogonal polynomial ensembles} \label{sec2}
\subsection{Notation}
The set of non-negative integers is denoted by $\mathbb{N}$. For $n \in \mathbb{N}$ and $a \in \mathbb{C}$, $(a)_n$ is the Pochhammer symbol:
\begin{align*}
(a)_0= 1 \quad ; \quad  (a)_n := a(a+1)...(a+n-1)=\frac{\Gamma(a+n)}{\Gamma(a)}, \quad n = 1,2,...
\end{align*}
where $\Gamma$ is the Euler Gamma function.\\

For $a_1,a_2,...,a_p \in \mathbb{C}$, $b_1,...,b_q \in \mathbb{C} \setminus (-\mathbb{N}) $ the hypergeometric function $\left._pF_q \right.(a_1,...,a_p ; b_1,...,b_q ; . )$ is defined by 
\begin{align} \label{defgauss}
\left._pF_q \right.(a_1,...,a_p ; b_1,...,b_q ; z) = \sum_{n=0}^{+\infty} \frac{(a_1)_n...(a_p)_n}{(b_1)_n...(b_q)_n} \frac{z^n}{n!} , \quad z \in \mathbb{C},
\end{align}
provided the series converges. The function $\left._2F_1 \right.$ is the Gauss hypergeometric function.
\\

If $\lambda=(\lambda_1 \geq \lambda_2 \geq ... ) \in \mathbb{Y}$ is a partition and $a \in \mathbb{C}$ is a complex number, we define the generelized Pochhammer symbol $(a)_\lambda$ by:
\begin{align*}
(a)_\lambda:=\prod_{i=1}^{l(\lambda)}(a-i+1)_{\lambda_i}=\prod_{(i,j) \in \lambda} (a-i+j)
\end{align*}
where we write 
\begin{align*}
(i,j) \in \lambda \Leftrightarrow 1 \leq j \leq \lambda_i.
\end{align*}
\subsection{Christoffel deformations of discrete orthogonal polynomials}
Let $\omega$ be a discrete measure on $\mathbb{R}$ with infinite support and finite moments of all orders:
\begin{align*}
\int_{\mathbb{R}}|x|^nd\omega(x) < +\infty \quad \text{for all $n  \in \mathbb{N}$}.
\end{align*}
We denote by $\{p_n\}_{n\in \mathbb{N}}$ a family of orthogonal polynomials with respect to $\omega$. The number $c_n$ stands for the leading coefficient of $p_n$, and its squared norm will be denoted by $h_n$:
\begin{align*}
p_n(x)&= c_nx^n + \text{terms of degree less than $n$}, \\
\int_{\mathbb{R}} p_n(x)p_m(x) d\omega(x) & = \begin{cases} 0 \quad \text{if $n \neq m$},\\
 h_n \quad \text{if $n=m$}
 \end{cases}.
\end{align*}
Let $k$ be a non-negative integer and let $u_1$,...,$u_k$ be pairwise distinct real numbers.
\begin{defi}The Christoffel deformation $\omega^k$ of $\omega$ of order $k$ at points $u_1$,...,$u_k$ is the discrete measure on $\mathbb{R}$ defined by:
\begin{align*}
d\omega^k(x) = \prod_{j=1}^k (x-u_j)^2d\omega(x)
\end{align*}
\end{defi}
Let $\{p_n^k \}_{n \in \mathbb{N}}$ be the family of monic orthogonal polynomials with respect to $\omega^k$. We have the following 
\begin{prop} \label{propchristoffel} The following explicit formula holds:
\begin{align*} 
p_n^k(x)=\left(\prod_{i=1}^k(x-u_i)^2.\delta_n^k.c_{n+2k}\right)^{-1}D_n^k(x)
\end{align*}
where
		\begin{align*}
		D_n^k(x)&= \begin{vmatrix}
		p_n(u_1) & p_{n+1}(u_1) & \dots & p_{n+2k}(u_1) \\			\vdots & \vdots & \ddots & \vdots \\
		p_n(u_k) & \dots & \dots &  p_{n+2k}(u_k) \\
		p'_n(u_1) & p'_{n+1}(u_1) & \dots & p'_{n+2k}(u_1) \\
		\vdots & \vdots & \ddots & \vdots \\
		p'_n(u_k) & \dots & \dots &  p'_{n+2k}(u_k) \\
		p_n(x) & p_{n+1}(x) & \dots & p_{n+2k}(x)
		\end{vmatrix} ; \\
		\\
		\delta_n^k &= \begin{vmatrix}
		p_n(u_1) & p_{n+1}(u_1) & \dots & p_{n+2k-1}(u_1) \\			\vdots & \vdots & \ddots & \vdots \\
		p_n(u_k) & \dots & \dots  &  p_{n+2k-1}(u_k) \\
		p'_n(u_1) & p'_{n+1}(u_1) & \dots & p'_{n+2k-1}(u_1) \\
		\vdots & \vdots & \ddots & \vdots \\
		p'_n(u_k) & p'_{n+1}(u_k) & \dots & p'_{n+2k-1}(u_k)
		\end{vmatrix}.
		\end{align*}
		Its squared $l^2(\omega^k)$-norm is given by
\begin{align*}
h_n^k=\frac{h_n\delta_{n+1}^k}{\delta_n^k.c_{n+2k}.c_n}
\end{align*}
\end{prop}
\begin{proof}
We slightly adapt the proof from \cite{ismail}, where it is assumed that the $u_i$ are outside the support of the measure $\omega$, to the case when the measure $\omega$ has infinite support, see also \cite{szego} and \cite{baikdeiftstrahov}. We first prove that $\delta_n^k \neq 0$. Let $\alpha_1, \dots, \alpha_{2k} \in \R$ be $2k$ real numbers such that the polynomial
\[P(x):= \sum_{i=1}^{2k} \alpha_i p_{n+i-1}(x) \]
and its derivative vanish for $x \in \{u_1, \dots , u_k\}$. Then there exists a polynomial $Q$ of degree at most or equal to $n-1$ such that
\[P(x)= Q(x) \prod_{i=1}^k (x-u_i)^2. \]
Since the polynomials $p_{n +i-1}$, $i=1,\dots,2k$ are orthogonal to any polynomial of degree less than or equal to $n-1$ with respect to $\omega$, we have
\begin{align} \label{eq:orthochristoffel1} \int_\R Q(x) P(x) d\omega(x) = \int_\R Q^2(x) d\omega^k(x) = 0.
\end{align}
Since the support of the measure $\omega$ is infinite, so is the support of the measure $\omega^k$, and thus Equation (\ref{eq:orthochristoffel1}) implies $Q=0$ and so $P=0$ i.e. $\alpha_1=\dots=\alpha_{2k}=0$. This implies that $\delta_n^k \neq 0$.

The fact that $\delta_n^k \neq 0$ tells us that $D_n^k(x)$ is a polynomial of degree $n+2k$. Since $D_n^k(x)$ and its derivative vanish for $x \in \{u_1,\dots,u_k \}$, the rational function $p_n^k(x)$ is in facts a monic polynomial of degree $n$.

Let $q(x)$ be a polynomial of degree strictly less than $n$. From the orthogonality of the polynomials $p_n$ with respect to $\omega$, we have
\begin{align*}
\int_\R p_n^k(x) q(x) d\omega^k(x) = \frac{1}{c_{n+2k}\delta_n^k} \int_\R D_n^k(x) q(x) d\omega(x) = 0,
\end{align*}
and the polynomials $p_n^k$ are thus orthogonal with respect to $\omega^k$. The squared norm $h_n^k$ is computed in \cite[Prop. 4.2.]{bufetovqiupickrell}.
\end{proof}
\subsection{Discrete orthogonal polynomial ensembles and determinantal point processes}
Let $E$ be an infinite countable discrete subset of $\mathbb{R}$. We first recall the definition of an orthogonal polynomial ensemble in our setting. Let $\omega$ be a discrete weight as in the previous subsection. We assume that its support lies inside $E$, and in all the sequel, we may identify a measure on $E$ with its density with respect to the counting measure on $E$. Let $N$ be a positive integer.
\begin{defi}An $N$-orthogonal polynomial ensemble is a probability measure $\mathbb{P}_N$ on $E^N$ given by the following formula:
\begin{equation} \label{orthpolenseq}
d\mathbb{P}_N(x_1,...x_N)= C_N \prod_{1 \leq i < j \leq N} (x_i-x_j)^2 \prod_{i=1}^Nd\omega(x_i)
\end{equation}
where $C_N$ is the normalization constant:
\begin{align*}
C_N^{-1} = \sum_{x_1, ..., x_N \in E} \prod_{1 \leq i < j \leq N} (x_i-x_j)^2 \prod_{i=1}^N\omega(x_i) 
\end{align*}
\end{defi}
We now recall the definition of a determinantal point process on $E$. We denote by $\text{Conf}(E)$ the space of configuration in $E$ which is in this context simply the set of all subsets of $E$,
\begin{align*}
\text{Conf}(E) := 2^E.
\end{align*}
\begin{defi} \label{defdpp}A \textit{point process} is a probability measure $\mathbb{P}$ on $\text{Conf}(E)$. It is said to be \textit{determinantal} if there exists a function
\begin{align*}
K : E \times E \rightarrow \mathbb{R}
\end{align*}such that for any $n \in \mathbb{N}$ and any $\{a_1,...,a_n \} \subset E$, one has 
\begin{align*}
\mathbb{P}( X \in \text{Conf}(E), \hspace{0.1cm}  \{a_1,...,a_n \} \subset X) = \det \left(K(a_i,a_j)\right)_{i,j=1}^n.
\end{align*}
The function $K$ is called the \textit{correlation kernel} of the determinantal point process.
\end{defi}
\begin{rem}The correlation kernel of a given determinantal point process is not unique. Indeed, if one has a determinantal point process with kernel $K$, for any non-vanishing function $f:E \rightarrow \mathbb{R}$ the function
\begin{align*}
\tilde{K}(x,y)=\frac{f(x)}{f(y)}K(x,y)
\end{align*}
can also serve as a kernel for the same determinantal point process.
\end{rem}
As a matter of fact, any orthogonal polynomial ensemble gives rise to a determinantal point process, what we resume in the following Proposition (see e.g. \cite{konig} and references therein). Note that since the measure (\ref{orthpolenseq}) is symmetric and does not charge $N$-tuples which contain two equal coordinates, it can and it will be identified with a measure on the space of configurations.
\begin{prop} \label{proporthpolens} An orthogonal polynomial ensemble given by (\ref{orthpolenseq}) is a determinantal point process on $E$. Its correlation kernel is given by the formula
\begin{align*}
K_N(x,y)&=\sqrt{\omega(x)\omega(y)}\sum_{n=0}^{N-1}\frac{p_n(x)p_n(y)}{h_n}\\
&= \begin{cases}\sqrt{\omega(x)\omega(y)}\frac{c_{N-1}}{h_{N-1} c_{N}}\frac{p_N(x)p_{N-1}(y)-p_N(y)p_{N-1}(x)}{x-y}, \quad \text{for $x \neq y$}, \\
\omega(x)\frac{c_{N-1}}{h_{N-1} c_{N}}\left(p_{N}'(x)p_{N-1}(x)-p_N(x)p_{N-1}'(x)\right), \quad \text{for $x=y$}. \end{cases}
\end{align*}
\end{prop}
The second equality is just the Christoffel-Darboux formula, while the last one follows from L'Hospital rule.
\begin{rem}Observe that $K_N$ is the kernel of the orthogonal projection from $l^2(E)$ (where $E$ is equipped with the counting measure) onto 
\begin{align*}\sqrt{\omega}\mathbb{R}_{N-1}[x] = \{x \mapsto\sqrt{\omega(x)}p(x),\hspace{0.1cm} p \text{ is a polynomial of degree at most $N-1$}\}.
\end{align*}
\end{rem}
The main object of our study is given by the following definition.
\begin{defi}[Christoffel deformations of orthogonal polynomial ensembles] \label{maindef}
Let $k \in \mathbb{N}$ be a non-negative integer. Given an orthogonal polynomial ensemble $\mathbb{P}_N$ with weight $\omega$, its \textit{Christoffel deformation} of order $k$, denoted by $\mathbb{P}_N^{k}$, is the orthogonal polynomial ensemble with weight $\omega^k$.
\end{defi}
The next Proposition gives a compact formula for the kernel of the Christoffel deformation which will be suitable for asymptotic analysis.
\begin{prop}\label{propkernelcharlier}The Christoffel deformation of order $k$, $\mathbb{P}_N^k$, is a determinantal point process with kernel $K_N^k$ given by:
\begin{align*}
K^{k}_N(x,y) =\frac{\sqrt{\omega(x)\omega(y)}}{\prod_{i=1}^k|(x-u_i)(y-u_i)|} \frac{c_{N-1}}{h_{N-1}c_{N+2k}}\frac{D_N^k(x)D_{N-1}^k(y)-D_N^k(y)D_{N-1}^k(x)}{(\delta_N^k)^2(x-y)}.
\end{align*}
\end{prop}
\begin{proof}This is a straightforward computation, using  Propositions \ref{proporthpolens} and \ref{propchristoffel}. Indeed, Proposition \ref{proporthpolens} states that $\mathbb{P}_N^k$ is a determinantal point process with kernel 
\begin{align*}
K_N^k(x,y)=\sqrt{\omega^k(x)\omega^k(y)}\frac{p_N^k(x)p_{N-1}^k(y)-p_N^k(y)p_{N-1}^k(x)}{h_{N-1}^k(x-y)}.
\end{align*}
By the definition of the measure $\omega^k$ and using Proposition \ref{propchristoffel} and the explicit expression of the polynomials $p_n^k$, this kernel can be written as 
\begin{multline*}
K_N^k(x,y)\\
=\frac{\sqrt{\omega(x)\omega(y)}}{\prod_{i=1}^k|(x-u_i)(y-u_i)|}\frac{1}{\delta_N^k\delta_{N-1}^k c_{N+2k}c_{N+2k-1}h_{N-1}^k}\frac{D_N^k(x)D_{N-1}^k(y)-D_N^k(y)D_{N-1}^k(x)}{x-y}.
\end{multline*}
The expression for the squared norm $h_{N-1}^k$ given in Proposition \ref{propchristoffel} leads to the desired result.
\end{proof}
\subsection{Palm measures of determinantal point processes} \label{secpalm}
We here give a probabilistic interpretation of Christoffel deformations of orthogonal polynomial ensembles in the case when the points $u_1,...,u_k$ belong to the support of the measure $\omega$. We first recall the basics on Palm distributions of a general point process $\P$ on a discrete space $E$, see for instance \cite{daley-verejones}, \cite{shiraitakahashi} or \cite{bufetovquasisymmetries} and references therein for a more general treatment.
\begin{defi}Let $\P$ be a point process on the discrete set $E$. For pairwise distinct points $u_1,...,u_k \in E$ such that $\P(u_1,...,u_k  \in X) >0$, the \textit{Palm measure} $\hat{\P}^{u_1,...,u_k}$ of $\P$ at $u_1,...,u_k$ is the point process $\P$ on $E$ conditionned to have particles at $u_1,...,u_k$,
\begin{align*}
\hat{\P}^{u_1,...,u_k} ( \mathcal{B} )= \frac{\P(\mathcal{B}\cap \{u_1,...,u_k  \in X \})}{\P(u_1,...,u_k  \in X) },
\end{align*}
for any borel set $ \mathcal{B} \subset \Conf(E)$. The \textit{reduced Palm measure} $\P^{u_1,...,u_k}$ is then defined as the pushforward of the Palm measure $\hat{\P}^{u_1,...,u_k}$ under the map
\begin{align*}
\Conf(E) & \rightarrow \Conf(E)\\
X & \mapsto X \setminus \{u_1,...,u_k \}.
\end{align*}
\end{defi}
In particular, if the point process $\P$ has $N + k$ particles, i.e. $\P( |X| = N+k) =1$, then its reduced Palm measure $\P^{u_1,\dots, u_k}$ has $N$ particles.\\

We also recall the following Proposition from \cite{shiraitakahashi}.
\begin{prop} \label{propcvpalm}
If $\mathbb{P}$ is a determinantal point process with kernel $K$ such that $\P(u_1,...,u_k  \in X) >0$, then, its reduced Palm measure at $u_1,...,u_k$ is a determinantal point process with kernel
\begin{align*}
K^{u_1,...,u_k}(x,y)= \left( \det(K(u_i,u_j)\right)_{1 \leq i,j \leq k}^{-1} \det \begin{pmatrix}
K(x,y) & K(x,u_1) &\dots & K(x,u_k) \\
K(u_1,y) &K(u_1,u_1) & \dots & K(u_1,u_k) \\
\vdots & & \ddots & \vdots \\
K(u_k,y) & \dots & \dots& K(u_k,u_k)
\end{pmatrix}.
\end{align*}
acting on $l^2(E )$.
\end{prop}
We now consider a measure $\omega$ as before, with an infinite support lying in a discrete countable subset $E$ of $\mathbb{R}$. For $n \in \mathbb{N}$, we denote by $\P_n$ the corresponding orthogonal polynomial ensemble. Let $u_1,...,u_k \in \text{supp}(\omega)$ be pairwise distinct, and, as before, denote by $\P_n^k$ the Christoffel deformation of $\P_n$ at $u_1,...,u_k$. The next Proposition formally establishes the fact that the study of Christoffel deformations of orthogonal polynomial ensembles contains the study of their Palm measures.
\begin{prop} \label{proppalm}
For any natural $N \in \mathbb{N}$, we have the equality between point processes
\begin{align*}
\P_{N+k}^{u_1,...,u_k}=\P_N^k.
\end{align*}
\end{prop}
\begin{proof}We know from \cite{shiraitakahashi} (see also \cite{bufetovquasisymmetries}) that, for a determinantal point process with a kernel that is a kernel of an orthogonal projection onto a closed subspace $L \subset l^2(E)$, its reduced Palm measure at points $u_1,...,u_k$ is a determinantal point process with the kernel of orthogonal projection onto
\begin{align*}
L^k:= \{f \in L, \hspace{0.1cm} f(u_1)=...=f(u_k)=0 \}.
\end{align*}
Specializing this fact to the subspace
\begin{align*}
L:= \{ \sqrt{\omega(x)}p(x), \hspace{0.1cm} p \in \mathbb{R}_{N+k-1}[x] \},
\end{align*} one obtains the statement of the Proposition. \end{proof}
\section{The Charlier ensemble and a deformation of the Plancherel measure} \label{sec3}
\subsection{On the Plancherel measure on partitions}
The Plancherel measure is a probability measure $\text{Pl}_n$ on the set of Young diagrams of fixed size $n$, $\mathbb{Y}_n=\{\lambda \in \mathbb{Y}, \hspace{0.1cm} |\lambda|=n \}$. It is given by 
\begin{align*}
\text{Pl}_n(\lambda)= \frac{\dim(\lambda)^2}{n!}, \quad \lambda \in \mathbb{Y}_n,
\end{align*}
where $\dim(\lambda)$ is the dimension of the irreducible representation of the $n$-th symmetric group parametrized by $\lambda$, or equivalentely, the number of stantard Young tableaux of shape $\lambda$ (see \cite{macdonald} or \cite{sagan}). It is well known that the poissonization of the Plancherel measure is a determinantal point process. More precisely, let $\alpha > 0$ be a positive parameter and form the poissonized Plancherel measure $\mathbb{P}_\alpha$ on $\mathbb{Y}$ by imposing $n$ to have the Poisson distribution with parameter $\alpha$:
\begin{align*}
\mathbb{P}_\alpha (\lambda) =e^{-\alpha} \alpha^{|\lambda|}\left(\frac{\dim(\lambda)}{|\lambda|!}\right)^2, \quad \lambda \in \mathbb{Y}.
\end{align*} For $\lambda \in \mathbb{Y}$, we set 
\begin{align*}
\mathfrak{S}^0(\lambda):=\{\lambda_i -i, \hspace{0.1cm} i=1,2,...\} \in \text{Conf}(\mathbb{Z}).
\end{align*}Then the pushforward of $\mathbb{P}_\alpha$ under the map $\mathfrak{S}^0$, which we call the \textit{discrete Bessel point process}, is a determinantal point process with the discrete Bessel kernel 
\begin{equation} \label{discretebesselkernel} K_{\alpha}(x,y)=\sqrt{\alpha} \frac{J_x(2\sqrt{\alpha})J_{y+1}(2\sqrt{\alpha})-J_y(2\sqrt{\alpha})J_{x+1}(2\sqrt{\alpha})}{x-y},
\end{equation}
where the Bessel functions $J_x(2\sqrt{\alpha}), \hspace{0.1cm} x \in \mathbb{Z}$ are defined via their generating Laurent series (\cite{andrews}):
\begin{equation} \label{defbessel}
e^{\sqrt{\alpha}(z-z^{-1})}=\sum_{x \in \mathbb{Z}} J_x(2\sqrt{\alpha}) z^x.
\end{equation}
This definition leads to the expression
\begin{align*}
J_x(2\sqrt{\alpha}) = \sum_{m=0}^{+ \infty}\frac{(-1)^m}{m!\Gamma(m+x+1)} \alpha^{m + x/2}
\end{align*}
that can be analytically continuated for all values of the indices $x \in \mathbb{C}$. We can thus introduce the derivatives of the Bessel functions with respect to their index 
\begin{align*}
\frac{d}{dx}J_x(2\sqrt{\alpha})=:L_x(2\sqrt{\alpha}).
\end{align*}
\subsection{Christoffel deformations of the Charlier ensemble and main result}We now define the Christoffel deformation of the Charlier ensemble. We know from \cite{johanssonplancherel} that this ensemble approximates the poissonized Plancherel measure, and the aim of this section is to generalize this fact by proving that Christoffel deformations of this ensemble admit a limit that is expressed through Bessel functions. Our discrete measure is 
\begin{align*}
\omega(x)=\omega_a(x)=e^{-a}\frac{a^x}{x!}, \quad x \in \mathbb{N},
\end{align*}
where $a>0$ is a parameter, and the $N$-th orthogonal polynomial ensemble that this measure defines is denoted by $\mathbb{P}_{N,a}$. The corresponding orthonormal polynomials are the Charlier polynomials (\cite{ismail}, \cite{andrews}):
\begin{align*}
C_n(x)=C_n(x;a)=\frac{a^{\frac{n}{2}}}{\sqrt{n!}}\left._2F_0\right.(-n,-x; ;-1/a) \hspace{0,1cm} , \quad n \in \mathbb{N}.
\end{align*}
As in the previous section, we define the Christoffel deformations of the Charlier ensemble. Let $k \in \mathbb{N}$ be a non-negative integer, and $\tilde{u}_i \in \mathbb{R}$ be distinct numbers.  Let $\alpha >0$ be a positive real number, and for each $N \in \mathbb{N}$, consider the Christoffel deformation of the $N$-th Charlier ensemble at points $u_1=\tilde{u}_1+N,...,u_k=\tilde{u}_k+N$, with parameter $a=\alpha/N$, according to definition \ref{maindef}. This determinantal point process is denoted by $\mathbb{P}_{N,\alpha/N}^k$, and its correlation kernel is denoted by $K_N^{\alpha/N,k}$. For an integer $p \in \mathbb{Z}$, we form the determinants:
\begin{align*}
A_{k,p}(x) &= \begin{vmatrix}
J_{\tilde{u}_1+p}(2\sqrt{\alpha}) &\dots & J_{\tilde{u}_1-2k+p}(2\sqrt{\alpha}) \\
\vdots& \ddots & \vdots \\
J_{\tilde{u}_k+p}(2\sqrt{\alpha}) & \dots& J_{\tilde{u}_k-2k+p}(2\sqrt{\alpha}) \\
L_{\tilde{u}_1+p}(2\sqrt{\alpha}) &\dots & L_{\tilde{u}_1-2k+p}(2\sqrt{\alpha}) \\
\vdots & \ddots & \vdots \\
L_{\tilde{u}_k+p}(2\sqrt{\alpha}) & \dots & L_{\tilde{u}_k-2k+p}(2\sqrt{\alpha}) \\
J_{x+p}(2\sqrt{\alpha}) & \dots & J_{x-2k+p}(2\sqrt{\alpha})
\end{vmatrix}, \quad x \in \mathbb{Z},\\
C_{k,p} &=\begin{vmatrix}
J_{\tilde{u}_1+p}(2\sqrt{\alpha}) & \dots & J_{\tilde{u}_1-2k+1+p}(2\sqrt{\alpha}) \\
\vdots & \ddots & \vdots \\
J_{\tilde{u}_k+p}(2\sqrt{\alpha}) &\dots& J_{\tilde{u}_k-2k+1+p}(2\sqrt{\alpha}) \\
L_{\tilde{u}_1+p}(2\sqrt{\alpha}) & \dots & L_{\tilde{u}_1-2k+1+p}(2\sqrt{\alpha}) \\
\vdots & \ddots & \vdots \\
L_{\tilde{u}_k+p}(2\sqrt{\alpha}) & \dots & L_{\tilde{u}_k-2k+1+p}(2\sqrt{\alpha})
\end{vmatrix}.
\end{align*}
The main result of this section is the following Theorem.
\begin{thm} \label{thm1}
In the regime described above, we have that, for all $x,y \in \mathbb{Z}$, $x \neq y$,
\begin{align*}
\lim_{N \rightarrow \infty} K_N^{\alpha/N,k}(x+N,y+N)&=\frac{\alpha^{\frac{2k+1}{2}}}{C_{k,0}^2 \prod_{i=1}^k|(x-\tilde{u}_i)(y-\tilde{u}_i)|} \frac{A_{k,0}(x)A_{k,1}(y)-A_{k,0}(y)A_{k,1}(x)}{x-y}.
\end{align*}
When $x=y$, L'Hospital rule entails:
\begin{align*}
\lim_{N \rightarrow \infty} K_N^{\alpha/N,k}(x+N,x+N)&=\frac{\alpha^{\frac{2k+1}{2}}}{C_{k,0}^2 \prod_{i=1}^k(x-\tilde{u}_i)^2}\left( A_{k,0}'(x)A_{k,1}(x)-A_{k,0}(x)A_{k,1}'(x)\right).
\end{align*}
\end{thm}
This implies that the Christoffel deformation of the Charlier ensemble $\mathbb{P}_{N,\alpha/N}^k$, shifted by $-N$, weakly converges, as a probability measure on $\text{Conf}(\mathbb{Z})$, to a determinantal point process governed by the limit correlation kernel of the Theorem. In particular, if we take $k=0$, one recovers Johansson's result (\cite{johanssonplancherel}, Theorem 1.2).
\subsection{Proof of Theorem \ref{thm1}}
We will use the following Lemma, which we prove at the end of this section. This Lemma gives simple asymptotic results on the Charlier polynomials involving Bessel functions and their derivatives. As far as we know, the asymptotics we give can not be found in the litterature in this simple form. However, one may consult \cite{ouwong} for a complete study of the asymptotics of Charlier polynomials.
\begin{lem}\label{lemme1}
For any $l \in \mathbb{Z}$ and any $u \in \mathbb{R}$, we have
\begin{align*}
\lim_{N \rightarrow \infty}\sqrt{\omega_{\alpha/N}(u+N)} C_{N+l}(u+N;\alpha/N)=J_{u-l}(2\sqrt{\alpha}),
\end{align*}
and
\begin{align*}
\lim_{N \rightarrow \infty} \sqrt{\omega_{\alpha/N}(u+N)} \left( C'_{N+l}(u+n;\alpha/N) - \log \left(N/\sqrt{\alpha} \right) C_{N+l}(u+N;\alpha/N) \right) = L_{u-l}(2\sqrt{\alpha}).
\end{align*}
\end{lem}
The proof of Theorem \ref{thm1} goes as follows. We use the formula from Proposition \ref{propkernelcharlier} for the kernel $K_N^{\alpha/N, k}$ and proceed to the analysis of each factor. First, a straightforward computation shows that:
\begin{equation}\label{limconstante}
\frac{c_{N-1}}{c_{N+2k}} \rightarrow \alpha^{\frac{2k+1}{2}}
\end{equation}
when $N \rightarrow +\infty$ and where $c_n$ is the leading coefficient of the Charlier polynomial $C_n$. We now analyse the ratio
\begin{align*}
\sqrt{\omega_{\alpha/N}(x+N)\omega_{\alpha/N}(y+N)}\frac{D_N^k(x+N)D_{N-1}^k(y+N)-D_N^k(y+N)D_{N-1}^k(x+N)}{\left(\delta_N^k\right)^2(x-y)}.
\end{align*}
To this end, we multiply both sides of this ratio by
\begin{align*}m:=\prod_{i=1}^{k}\omega_{\alpha/N}(\tilde{u}_{i}+N)^2,
\end{align*} and in each determinant $D_N^k(.)$, $D_{N-1}^k(.)$ and $\delta_N^k$, we perform the following operations on rows $R_l$:
\begin{align*}
R_{k+i} \leftarrow R_{k+i}-\log \left( N/\alpha \right) R_{i}
\end{align*} for $i=1,...,k$. Let $M$ be the diagonal $2k \times 2k$ matrix:
\begin{align*}M=& \text{diag} \left( \sqrt{\omega_{\alpha/N}(u_1)},...,\sqrt{\omega_{\alpha/N}(u_k)},\sqrt{\omega_{\alpha/N}(u_1)}, ...,\sqrt{\omega_{\alpha/N}(u_k)} \right).
\end{align*}
Then
\begin{align*}
m=\det(M^2)=\det (M)^2.
\end{align*}Let us denote by $\hat{\delta_N^k}$ the matrix such that $\det\left( \hat{\delta_N^k} \right)=\delta_N^k$ with changed rows as described above and use the same notation for the other determinants. We have:
\begin{align*}
m\left( \delta_N^k \right)^{2}=\det \left(M\hat{\delta_N^k} \right)^2
\end{align*}
and each term of this matrix is precisely of the type we can analyse by Lemma \ref{lemme1}. In the same fashion, if we denote by $M_z$ the $(2k+1) \times (2k+1)$ matrix 
\begin{align*}
M_z=\text{diag}\left(M,\sqrt{\omega_{\alpha/N}(z+N)}\right),
\end{align*}
we have that 
\begin{multline*}
\sqrt{\omega_{\alpha/N}(x+N)\omega_{\alpha/N}(y+N)}\cdot m \cdot \left(D_{N}^k(x+N)D_{N-1}^k(y+N)-D_{N}^k(x+N)D_{N-1}^k(y+N) \right) \\
=\det (M_x) \det (M_y) \cdot \left(D_{N}^k(x+N)D_{N-1}^k(y+N)-D_{N}^k(x+N)D_{N-1}^k(y+N) \right) \\
= \det \left( M_x \hat{D}_N^k(x+N)\right) \det \left( M_y \hat{D}_{N-1}^k(y+N)\right)
\\
-\det \left(M_y \hat{D}_N^k(y+N) \right) \det \left( M_x \hat{D}_{N-1}^k(x+N) \right),
\end{multline*}
and again, we know the asymptotics for the components of each matrix thanks to Lemma \ref{lemme1}. Using the continuity of the determinant, we have established that, for all $x \neq y$ 
\begin{align*}
\lim_{N \rightarrow +\infty} \sqrt{\omega_{\alpha/N}(x+N)\omega_{\alpha/N}(y+N)}&\frac{D_N^k(x+N)D_{N-1}^k(y+N)-D_N^k(y+N)D_{N-1}^k(x+N)}{\left(\delta_N^k\right)^k(x-y)}\\
&=\frac{A_{k,0}(x)A_{k,1}(y)-A_{k,0}(y)A_{k,1}(x)}{C_{k,0}(x-y)}.
\end{align*}

The case $x=y$ goes as follows. For a $(2k+1)\times(2k+1)$ determinant $A$, and for $i=0,...,2k$ we denote by $[A]^i$ its $\left((i+1),(2k+1)\right)$ cofactor, i.e. the same determinant with the $(i+1)$-th column and the last line removed, multiplied by $(-1)^{i+1+2k+1}=(-1)^i$. We have 
\begin{multline} \label{eqsym}
D_{N}^k(x+N)'D_{N-1}^k(x+N)-D_{N}^k(x+N)D_{N-1}^k(x+N)'\\
=\left(\sum_{i=0}^{2k}\left[D_N^k(x+N)\right]^iC_{N+i}'(x+N) \right)\left(\sum_{i=0}^{2k}\left[D_{N-1}^k(x+N)\right]^iC_{N-1+i}(x+N) \right)\\
-\left(\sum_{i=0}^{2k}\left[D_N^k(x+N)\right]^iC_{N+i}(x+N) \right)\left(\sum_{i=0}^{2k}\left[D_{N-1}^k(x+N)\right]^iC_{N-1+i}'(x+N) \right)\\
=\sum_{i,j=0}^{2k}\left[D_N^k(x+N)\right]^i\left[D_{N-1}^k(x+N)\right]^j
 \left(C_{N+i}'(x+N)C_{N-1+j}(x+N)-C_{N+i}(x+N)C_{N-1+j}'(x+N) \right).
\end{multline}
Remark that the analysis of the cofactors has been done, and observe now that we have nice symmetries in this expression. Indeed, for any $i,j=0,...,2k$, second part of Lemma \ref{lemme1} entails:
\begin{align*}
C_{N+i}&'(x+N)C_{N-1+j}(x+N)-C_{N+i}(x+N)C_{N-1+j}'(x+N)\\
&=\omega_{\alpha/N}(x+N)^{-1/2}\\
&\times\left\{ \left(L_{x-i}(2\sqrt{\alpha})+\log\left(N/\sqrt{\alpha}\right)C_{N+i}(x+N)+o(1)\right)C_{N-1+j}(x+N)\right.\\
&\left.-C_{N+i}(x+N)\left(L_{x-1+j}(2\sqrt{\alpha})+\log\left(N/\sqrt{\alpha}\right)C_{N-1+j}(x+N)+o(1)\right)\right\}.
\end{align*}
We see that the term involving the factor $\log(N/\sqrt{\alpha})$ vanishes. Multiplying (\ref{eqsym}) by $\omega_{\alpha/N}(x+N)$ and applying the first part of Lemma \ref{lemme1}, we are left with:
\begin{align*}
\omega_{\alpha/N}(x+N)&\left(D_{N}^k(x+N)'D_{N-1}^k(x+N)-D_{N}^k(x+N)D_{N-1}^k(x+N)'\right)\\
&=\sum_{i,j=0}^{2k}\left[D_N^k(x+N)\right]^i\left[D_{N-1}^k(x+N)\right]^j\\
&\times\left(L_{-x-i}(2\sqrt{\alpha})J_{-x+1-j}(2\sqrt{\alpha})-J_{-x-i}(2\sqrt{\alpha})L_{-x+1-j}(2\sqrt{\alpha}) + o(1)\right).
\end{align*}
This sum can be factorized, and we obtain that:
\begin{align*}
\omega_{\alpha/N}(x+N)& \cdot m \cdot \left(D_{N}^k(x+N)'D_{N-1}^k(x+N)-D_{N}^k(x+N)D_{N-1}^k(x+N)'\right)
\end{align*}
tends to
\begin{align*}
 A_{k,0}'(x)A_{k,1}(x)-A_{k,0}(x)A_{k,1}'(x)
\end{align*}
as $N$ tends to $+\infty$, which is the desired result. Recalling (\ref{limconstante}), the Theorem is proved. $\square$\\
\\
\begin{proof}[Proof of Lemma \ref{lemme1}]
Following \cite{johanssonplancherel}, we will use the integral representation for Charlier polynomials. Since the generating series for the normalized Charlier polynomials is:
\begin{align*}
\sum_{m=0}^{\infty} \left( \frac{\alpha}{N} \right)^{m/2}\frac{1}{\sqrt{m!}}C_m(x;\alpha/N)w^m=e^{-\alpha w / N}\left(1+w \right)^x,
\end{align*}
we have by Cauchy's formula, for $r>1$:
\begin{multline} \label{form:intrepcharlier}
C_m(x;\alpha/N)= \sqrt{m!} \left( \frac{N}{\alpha} \right)^{m/2}  \left( \frac{1}{2i\pi} \int_{|w|=r} e^{-\alpha w/N} (1+w)^{x} w^{-m} \frac{dw}{w} \right. \\
 \left. - (-1)^m \frac{\sin \pi x}{\pi} \int_{1}^{r} e^{\alpha u /N} (u-1)^{x} u^{-m} \frac{du}{u} \right).
\end{multline}
The second integral in (\ref{form:intrepcharlier}) arises as we take the principal branch of logarithm in order to define $(1+w)^x$ for $w \in \C \setminus (-\infty, -1]$, $x \in \R$. Performing the change of variables $w=Nz/\sqrt{\alpha}$, $u=Ns/\sqrt{\alpha}$, one obtains, for $r > \sqrt{\alpha}/N$:
\begin{align*}
C_{m}(x;\alpha/N)&=\sqrt{\frac{m!}{N^m}} \left( \frac{1}{2\pi}\int_{-\pi}^{\pi}e^{-\sqrt{\alpha}re^{i\theta}} \left( 1+ \frac{Nre^{i\theta }}{\sqrt{\alpha}} \right)^x \frac{d\theta}{(re^{i\theta})^m} \right. \\
&\left. -(-1)^m \frac{\sin \pi x}{\pi} \int_{\sqrt{\alpha} / N}^{r} e^{\sqrt{\alpha}s} \left( \frac{Ns}{\sqrt{\alpha}}-1 \right)^x s^{-m} \frac{ds}{s} \right).
\end{align*}
This expression can be suitably factorized:
\begin{equation}\label{eq2}
\begin{split}
C_{m}(x;\alpha/N)&=e^{-\sqrt{\alpha}} \sqrt{\frac{m!}{N^m}} \left( 1+ \frac{N}{\sqrt{\alpha}} \right)^{x}\\
&\times \left( \frac{1}{2\pi} \int_{-\pi}^{\pi}e^{\sqrt{\alpha}(1-re^{i\theta})} \left(  \frac{\sqrt{\alpha}+Nre^{i\theta }}{\sqrt{\alpha}+N} \right)^x \frac{d\theta}{(re^{i\theta})^m} \right.\\
&\left. -(-1)^m \frac{\sin \pi x}{\pi} \int_{\sqrt{\alpha} / N}^{r} e^{\sqrt{\alpha}(s+1)} \left( \frac{Ns-\sqrt{\alpha}}{\sqrt{\alpha}+N} \right)^x s^{-m} \frac{ds}{s} \right).
\end{split}
\end{equation}
Take $m=N+l$ and denote by $I_N^l(x)$ the sum of the integrals
\begin{equation} \label{eqintcharlier} \begin{split}
I_N^l(x)&:= \frac{1}{2\pi} \int_{-\pi}^{\pi}e^{\sqrt{\alpha}(1-re^{i\theta})} \left(  \frac{\sqrt{\alpha}+Nre^{i\theta }}{\sqrt{\alpha}+N} \right)^x \frac{d\theta}{(re^{i\theta})^{N+l}} \\
&-(-1)^{N+l} \frac{\sin \pi x}{\pi} \int_{\sqrt{\alpha} / N}^{r} e^{\sqrt{\alpha}(s+1)} \left( \frac{Ns-\sqrt{\alpha}}{\sqrt{\alpha}+N} \right)^x s^{-N-l} \frac{ds}{s}
\end{split},
\end{equation}
for any $l \in \mathbb{Z}$, and define
\begin{align*}
A_{\alpha}^N(x)=e^{-2\sqrt{\alpha}}\frac{N!}{N^N}\left(1+\frac{N}{\sqrt{\alpha}} \right)^{2x} \omega_{\alpha/N}(x).
\end{align*}
We have:
\begin{equation}\label{eq1}
\sqrt{\omega_{\alpha/N}(x)}C_{N+l}(x;\alpha/N)=\sqrt{\frac{(N+1)...(N+l)}{N^l}A_{\alpha}^N(x)}I_N^l(x).
\end{equation} 
Recall the integral representation for the Bessel function, which follows from (\ref{defbessel}):
\begin{equation} \label{eqintbessel}
J_{x}(2\sqrt{\alpha})= \frac{1}{2\pi}\int_{-\pi}^{\pi}e^{\sqrt{\alpha}(\frac{1}{r}e^{-i\theta}-re^{i\theta})}(re^{i\theta})^xd\theta-\frac{\sin\pi x}{\pi}\int_{0}^{r}e^{\sqrt{\alpha}(-1/s+s)}s^{x}\frac{ds}{s}.
\end{equation}
We first show that the integrals (\ref{eqintcharlier}) converges to the integrals (\ref{eqintbessel}), with a speed of order at least $1/N$, showing that the integrands converge. For the first integrand, we have by definition:
\begin{align*}
e^{\sqrt{\alpha}(1-z)}&\left( \frac{\sqrt{\alpha}+Nz}{\sqrt{\alpha}+N} \right)^{x+N}z^{-N-l} \\
&= \exp \left( \sqrt{\alpha}(1-z) -(N+l)\log z +(x+N)(\log(\sqrt{\alpha} +Nz)-\log(\sqrt{\alpha} +N)) \right).
\intertext{Factorizing inside the logarithms leads to:}
e^{\sqrt{\alpha}(1-z)}&\left( \frac{\sqrt{\alpha}+Nz}{\sqrt{\alpha}+N} \right)^{x+N}z^{-N-l} \\
&=\exp \left( \sqrt{\alpha}(1-z) -(l-x)\log z +(x+N)\left(\log \left(1+\frac{ \sqrt{\alpha}}{Nz} \right) - \log \left( 1+ \frac{\sqrt{\alpha}}{N} \right) \right) \right) \\
&=\exp \left( \sqrt{\alpha}(1/z-z)+O(1/N) \right)z^{x-l},
\end{align*}
where the last line follows from straightforward simplifications and first order expansions of the logarithms near $1$. In the same way, we have for the second integrand:
\begin{align*}
e^{\sqrt{\alpha}(s+1)} \left( \frac{Ns-\sqrt{\alpha}}{\sqrt{\alpha}+N} \right)^{x+N} s^{-N-l}=\exp \left(\sqrt{\alpha}(-1/s+s)+O(1/N)\right)s^{x-l}.
\end{align*}
We have shown that there exists a constant $K=K(\alpha,x,r,l)$ such that:
\begin{align*}
\left|e^{\sqrt{\alpha}(1-z)}\left( \frac{\sqrt{\alpha}+Nz}{\sqrt{\alpha}+N} \right)^{x+N}z^{-N-l}- \exp \left( \sqrt{\alpha}(1/z-z) \right)z^{x-l}\right| \leq \frac{K}{N}, \\
\intertext{and}
\left|e^{\sqrt{\alpha}(s+1)} \left( \frac{Ns-\sqrt{\alpha}}{\sqrt{\alpha}+N} \right)^{x+N} s^{-N-l}- e^{\sqrt{\alpha}(-1/s+s)}s^{x-l}\right| \leq \frac{K}{N},
\end{align*}
for all involved $z$ and $s$. This implies:
\begin{align*}
 | I_N^l(x+N)-J_{x-l}(2\sqrt{\alpha} )| \leq \frac{K}{N},
\end{align*}
as announced. Since as $N \rightarrow +\infty$, we have that:
\begin{align*}
A_{\alpha}^N (x+N) \rightarrow 1,
\end{align*}
and recalling (\ref{eq1}), we have established the first part of the Lemma. \\

We now move to the proof of the second statement of the Lemma. Since the techniques are very similar, we will give less details in our computations. Differentating expression (\ref{eq2}), we obtain:
\begin{align*}
C_m'(x;\alpha/N)&=\log \left( 1+ \frac{N}{\sqrt{\alpha}} \right) C_m(x;\alpha/N)\\
&+ e^{-\sqrt{\alpha}} \sqrt{\frac{m!}{N^m}} \left( 1+ \frac{N}{\sqrt{\alpha}} \right)^{x} \frac{d}{dx} \left( I^{m-N}_{N} (x) \right).
\end{align*} 
We thus only need to proof that the derivative of the integrals converges to the derivative of the Bessel function. Derivations under the sign of the integrals give
\begin{align*}
\frac{d}{dx} \left( I^{l}_{N} (x)\right)&=\frac{1}{2\pi} \int_{-\pi}^{\pi}\log \left( \frac{\sqrt{\alpha}+Nre^{i\theta }}{\sqrt{\alpha}+N} \right)e^{\sqrt{\alpha}(1-re^{i\theta})} \left(  \frac{\sqrt{\alpha}+Nre^{i\theta }}{\sqrt{\alpha}+N} \right)^x \frac{d\theta}{(re^{i\theta})^{N+l}} \\
&-(-1)^{N+l} \frac{\sin \pi x}{\pi} \int_{\sqrt{\alpha} / N}^{r} \log \left( \frac{Ns-\sqrt{\alpha}}{\sqrt{\alpha}+N} \right) e^{\sqrt{\alpha}(s+1)} \left( \frac{Ns-\sqrt{\alpha}}{\sqrt{\alpha}+N} \right)^x s^{-N-l} \frac{ds}{s} \\
&-(-1)^{N+l} \frac{\cos \pi x}{\pi} \int_{\sqrt{\alpha} / N}^{r} e^{\sqrt{\alpha}(s+1)} \left( \frac{Ns-\sqrt{\alpha}}{\sqrt{\alpha}+N} \right)^x s^{-N-l} \frac{ds}{s},
\end{align*} and
\begin{align*}
L_{x}(2\sqrt{\alpha})&= \frac{1}{2\pi}\int_{-\pi}^{\pi}\log(re^{i\theta})e^{\sqrt{\alpha}(\frac{1}{r}e^{-i\theta}-re^{i\theta})}(re^{i\theta})^xd\theta\\
&-\frac{\sin\pi x}{\pi}\int_{0}^{r}\log(s)e^{\sqrt{\alpha}(-1/s+s)}s^{x}\frac{ds}{s} \\
&-\frac{\cos\pi x}{\pi}\int_{0}^{r}e^{\sqrt{\alpha}(-1/s+s)}s^{x}\frac{ds}{s}.
\end{align*} Analogous computations lead to:
\begin{align*}
\lim_{N \rightarrow +\infty} \frac{d}{dx} \left( I^{l}_{N} (x+N)\right) = L_{x-l}(2\sqrt{\alpha}).
\end{align*}
This completes the proof of the Lemma.
\end{proof}
\subsection{Palm measures of the discrete Bessel point process}
In regards of subsection \ref{secpalm} and Proposition \ref{proppalm}, we can now give a description of the Palm measures of the discrete Bessel point process $\mathfrak{S}^0_*(\mathbb{P}_\alpha)$, the determinantal point process with kernel $K_\alpha$ given by (\ref{discretebesselkernel}).
\begin{thm} \label{thmpalmbessel}Let $\tilde{u}_1,...,\tilde{u}_k \in \mathbb{Z}$ be pairwise distinct. For $\alpha >0$, let $\left( \mathfrak{S}^0_* (\mathbb{P}_\alpha) \right)^{k}$ be the reduced Palm measure of the discrete Bessel point process at $\tilde{u}_1,...,\tilde{u}_k$. Then, the process $\left( \mathfrak{S}^0_* (\mathbb{P}_\alpha) \right)^{k}$ is a determinantal process with kernel
\begin{align*}
K_\alpha^{\tilde{u}_1,...,\tilde{u}_k}(x,y):=\frac{\alpha^{\frac{2k+1}{2}}}{C_{k,k}^2 \prod_{i=1}^k|(x-\tilde{u}_i)(y-\tilde{u}_i)|} \frac{A_{k,k}(x)A_{k,1+k}(y)-A_{k,k}(y)A_{k,1+k}(x)}{x-y}.
\end{align*}
\end{thm}
\begin{rem} \label{remtasep}From \cite{borodingorin}, Theorem 4.1, the kernel $K_\alpha^{\tilde{u}_1,...,\tilde{u}_k}$ serves as the correlation kernel of the TASEP with initial step condition at time $t=\alpha$ under the condition that the position $\tilde{u}_1,...,\tilde{u}_k$ are occupied.
\end{rem}
\begin{proof}[Proof of Theorem \ref{thmpalmbessel}]By Propositions \ref{proppalm} and \ref{propchristoffel}, the reduced Palm measure at $\tilde{u}_1+N,...,\tilde{u}_k+N$ of the $N$-th Charlier ensemble with parameter $\alpha/N$, shifted by $-N$, is a determinantal point process with kernel $K_{N-k}^{\alpha/N,k}(x+N,y+N)$. The case $k=0$ of Theorem \ref{thm1} together with Proposition \ref{propcvpalm} imply that the reduced Palm measure of the $N$-th Charlier ensemble at $\tilde{u}_1+N,...,\tilde{u}_k+N$ shifted by $-N$ converges to the reduced Palm measure of the discrete Bessel point process at points $\tilde{u}_1,...,\tilde{u}_k$. On the other hand, slight modifications of the proof of Theorem \ref{thm1} (precisely, we just replace $N$ by $N-k$) lead to
\begin{align*}
\lim_{N \rightarrow +\infty} K_{N-k}^{\alpha/N,k}(x+N,y+N)= K_{\alpha}^{\tilde{u}_1,...,\tilde{u}_k}(x,y).
\end{align*}
We deduce from this fact that the latter kernel $K_{\alpha}^{\tilde{u}_1,...,\tilde{u}_k}$ serves as a kernel for the point process $\left( \mathfrak{S}^0_*(\mathbb{P}_\alpha)\right)^{k}$.
\end{proof}
\section{The $z$-measures and a deformation of the gamma process} \label{sec4}
\subsection{Christoffel deformations of the $z$-measures}
The $z$-measures on partitions were introduced and extensively studied by Borodin and Olshanski in a series of papers, see e.g. \cite{bodist}, \cite{bogamma}, \cite{bomeixner} and \cite{bomarkov} and references therein. We first briefly explain their construcion. The $z$-measures are defined on the set $\mathbb{Y}$ of all partitions by:
\begin{align*}
M_{z,z',\xi}(\lambda)=(1-\xi)^{zz'}\xi^{|\lambda|}(z)_\lambda(z')_\lambda \left(\frac{\dim(\lambda)}{|\lambda|!}\right)^2.
\end{align*}
They depend on three complex parameters, $z$, $z'$ and $\xi$ and they are in general complex measures. The following choices ensure that it is a probability measure (see e.g. \cite{bomeixner} Proposition 1.2 for a proof). First of all, the parameter $\xi$ is real and lies in the open unit interval $(0,1)$. Next, we distinguish three cases:
\begin{enumerate}
	\item[$\bullet$][Principal series] The parameters $z$ and $z'$ are conjugate to each other.
	\item[$\bullet$][Complementary series] The parameters $z$ and $z'$ are both real and belong to the same interval $(m,m+1)$ where $m \in \mathbb{Z}$.
	\item[$\bullet$][Degenerate series] One of the parameter, say $z$, is a non-zero integer, while the other has the same sign and satisfy $ |z'|>|z|-1$.
\end{enumerate}
If $z$ and $z'$ satisfy one of these three conditions, we say that they are \textit{admissible} parameters.
The Christoffel defomations of the $z$-measure are then defined as follows:
\begin{defi}Let $z$, $z'$ be admissible parameters and $M_{z,z',\xi}$ the corresponding $z$-measure. Let $u_1$,...,$u_k$, $v_1$,...,$v_k \in\mathbb{C}$ be $2k$ complex numbers which satisfy:
\begin{align*}
\prod_{i=1}^k(x-u_i+z)(x-v_i+z') > 0
\end{align*}
for any $x \in \mathbb{Z}$. We define the Christoffel deformation of $M_{z,z',\xi}$ by:
\begin{align} \label{defzmeas}
M_{z,z',\xi}^k(\lambda)=C\prod_{j=1}^{l(\lambda)}\prod_{i=1}^k(\lambda_j-j-u_i+z)(\lambda_j-j-v_i+z')M_{z,z'\xi}(\lambda),
\end{align}
where $C$ is the normalisation constant chosen such that $M_{z,z',\xi}^k$ is a probability measure.
\end{defi}
We first have to focus on the degenerate case. For simplicity, we assume that $z=N$ is a positive integer and $z' > z-1$. The $z$-measure $M_{z,z'\xi}$ is then supported on the subset of Young diagrams of length less than or equal to $N$, and we map $\mathbb{Y}$ on $\text{Conf}(\mathbb{N}) \subset \text{Conf}(\mathbb{Z})$ via
\begin{align*}
\lambda \mapsto \mathfrak{S}_N(\lambda):=\{ \lambda_i-i+N, \hspace{0.1cm} i=1,...,N \}.
\end{align*}
The pushforward of the $z$-measure $M_{z,z'\xi}$ under $\mathfrak{S}_N$ is then the $N$-th Meixner ensemble (see e.g. \cite{bomeixner}, Proposition 1.4), i.e. the $N$-th orthogonal ensemble with weight
\begin{align*}
\omega_{\beta,\xi}(x)=\frac{(\beta)_x \xi^x}{x!},
\end{align*}
where $\beta = z'-z+1>0$. The corresponding orthogonal polynomials are the Meixner polynomials defined by:
\begin{align*}
M_n(x;\beta,\xi)= \left._2F_1\right.(-n,-x;\beta;1-1/\xi) , \hspace{0.1cm} n \in \mathbb{N}.
\end{align*} Their leading coefficient $c_n$ and their squared norm $h_n$ are given by
\begin{align*}c_n= \frac{\left(1-\frac{1}{\xi}\right)^n}{n!(\beta)_n} \hspace{0.1cm} , \quad h_n = \frac{n!}{\xi^n (1-\xi)^\beta (\beta)_n}.
\end{align*}
We now choose $k$ real numbers $u_1,...,u_k \in \mathbb{R}\setminus \mathbb{N}$ which are not integers, and set $v_i=u_i+\beta-1$. We thus have:
\begin{align*}
\prod_{j=1}^{l(\lambda)}\prod_{i=1}^k(\lambda_j-j-u_i+z)(\lambda_j-j-v_i+z')=\prod_{x \in \mathfrak{S}_N(\lambda)}\prod_{j=1}^k(x-u_i)^2
\end{align*}
from which we deduce that the pushforward of $M_{z,z',\xi}^k$ under $\mathfrak{S}_N$ is the Christoffel deformation of the $N$-th Meixner ensemble at points $u_1,...,u_k$ as defined in definition \ref{maindef}. As a consequence, it defines a determinantal point process with an explicit kernel given by Proposition \ref{propkernelcharlier}. We now prepare a similar statement for other choices of admissible parameters $z$ and $z'$. For this aim, we will reproduce the analytic continuation argument developed in \cite{bomeixner} and \cite{bomarkov}, although the simplifications we need are less obvious. We will note: $\mathbb{Z}'=\mathbb{Z}+1/2$ the set of proper half integers. For $a \in \mathbb{Z}'$, $x \in \mathbb{Z}'$ and $z$ and $z'$ admissible parameters \textit{that are not} in the degenerate case, we define:
\begin{equation}\label{defpsi} \begin{split}
\psi_a(x;z,&z',\xi):= \left( \frac{\Gamma(x+z+1/2) \Gamma(x+z'+1/2)}{\Gamma(z-a+1/2) \Gamma(z'-a+1/2)} \right)^{1/2} \xi^{1/2(x+a)} (1-\xi)^{1/2(z+z')-a} \\
&\times \frac{\left._2F_1\right.\left(-z+a+1/2,-z'+a+1/2;x+a+1;\frac{\xi}{\xi-1}\right)}{\Gamma(x+a+1)}.
\end{split}
\end{equation}
Observe that these functions are symmetric in $z,z'$. We now recall the following integral representation, which is due to Borodin-Olshanski (see \cite{bomeixner}, Proposition 2.3):
\begin{equation}\label{intreppsi} \begin{split}
\psi_a(x;z,z',\xi)&= \left( \frac{\Gamma(x+z+1/2)\Gamma(x+z'+1/2)}{\Gamma(z-a+1/2)\Gamma(z'-a+1/2)}\right)^{1/2} \frac{\Gamma(z'-a+1/2)}{\Gamma(z'+x+1/2)}(1-\xi)^{\frac{z'-z+1}{2}}\\
&\times\frac{1}{2i\pi} \oint (1- \sqrt{\xi}\omega)^{-z'+a-1/2} \left(1- \frac{\sqrt{\xi}}{\omega}\right)^{z-a-1/2}\omega^{-x-a}\frac{d\omega}{\omega}.
\end{split}
\end{equation}
Observe that it is now not obvious that such an expression is symmetric in $z,z'$, i.e. none of its factor is symmetric in $z,z'$, but since expression (\ref{defpsi}) is, it has to be symmetric in $z,z'$. We will also need the following Proposition, relating the functions $\psi_a$ to Meixner polynomials, which comes from Borodin and Olshanski (\cite{bomeixner}, Proposition 2.8).
\begin{prop} \label{propbopsimeixner}
Suppose $z=N \in \mathbb{N}$ and $z'=N+\beta-1$. Suppose that:
\begin{align*}
\tilde{x}:=x+N-1/2 \in \mathbb{N}, \quad \quad n:=N-a-1/2 \in \mathbb{N},
\end{align*}
are non-negative integers. Then, $\psi_a(x;z,z',\xi)$ is well defined and we have:
\begin{align*}
\psi_a(x;z,z',\xi) =\sqrt{\frac{\omega_{\beta,\xi}(\tilde{x})}{h_n}} M_n(\tilde{x};\beta,\xi).
\end{align*}
\end{prop}
In what follows, we will sometimes omit the dependance on the parameters $z$, $z'$ and $\xi$ when the context is clear. The first Theorem of this section states that, for any pair of admissible parameters $(z,z')$, the modified $z$-measure given by formula (\ref{defzmeas}) gives rise to a determinantal point process with an explicit kernel expressed through the functions $\psi_a$ and their derivatives. 
\begin{thm} \label{thmzmeas}For any admissible parameters $z$ and $z'$ which are not in the degenerate case, and any $u_i \in \mathbb{R}$, $v_i=u_i+z'-z-1/2$, the pushforward of the measure $M_{z,z',\xi}^k$ under the map
\begin{align*}
\mathfrak{S} : \lambda \mapsto \{\lambda_i-i+1/2\} \in \text{Conf}(\mathbb{Z}'),
\end{align*}
is a determinantal point process on $\mathbb{Z}'$. Its correlation kernel is given by:
\begin{equation} \label{kernelmodifiedzmeas}
K_{z,z',\xi}^k(x,y):=\frac{C_{z,z',k}}{\prod_{i=1}^{k}|(x-u_i)(y-u_i)|}\frac{B_{k,0}(x;z,z',\xi)B_{k,1}(y;z,z',\xi)-B_{k,0}(y;z,z',\xi)B_{k,1}(x;z,z',\xi)}{D_k(z,z',\xi)^2(x-y)},
\end{equation}
with
\begin{align*}
B_{k,p}(x)&= \begin{vmatrix}
 \psi_{-1/2+p}(u_1) & \psi_{-3/2+p}(u_1) & \dots & \psi_{-1/2-2k+p}(u_1) \\
 \vdots &   & \ddots & \vdots \\
 \psi_{-1/2+p}(u_k) & \dots & \dots & \psi_{-1/2-2k+p}(u_k) \\
 \psi'_{-1/2+p}(u_1) & \dots & \dots & \psi'_{-1/2-2k+p}(u_1) \\
 \vdots  & \ddots &  & \vdots \\
 \psi'_{-1/2+p}(u_k) & \dots & \dots & \psi'_{-1/2-2k+p}(u_k) \\
 \psi_{-1/2+p}(x) & \dots & \dots & \psi_{-1/2-2k+p}(x)
\end{vmatrix}, \quad p=0,1 \hspace{0.1cm}, \\
\\
D_k&=\begin{vmatrix}
 \psi_{-1/2}(u_1) & \psi_{-3/2}(u_1) &\dots& \psi_{1/2-2k}(u_1) \\
  \vdots & & \ddots & \vdots \\
 \psi_{-1/2}(u_k) & \dots & \dots & \psi_{1/2-2k}(u_k) \\
 \psi'_{-1/2}(u_1) & \dots & \dots & \psi'_{1/2-2k}(u_1) \\
 \vdots & \ddots &  & \dots \\
 \psi'_{-1/2}(u_k) & \dots & \dots & \psi'_{1/2-2k}(u_k) \\
\end{vmatrix}, \\
\intertext{and}
C_{z,z',k}&=\left(\frac{\xi}{\xi-1}\right)^{1+2k} \frac{\Gamma(z+2k+1)}{\Gamma(z)}\frac{\Gamma(z'+2k)}{\Gamma(z'-1)}.
\end{align*}
\end{thm}
\begin{proof}The proof requires some steps. We first establish in step 1 that the statement of the Theorem holds for the degenerate series case. This fact is precisely formulated in Lemma \ref{lemmestep1} below. In step 2, we prove that quantities of interest involving the modified $z$-measures have an analytic  in $\xi$ (and thus in $\sqrt{\xi}$), with polynomial coefficients in $(z,z')$, for all values of the admissible parameters. In step 3, we prove that a modification of the kernel also have an analytic expansion in $\sqrt{\xi}$ with polynomials coefficients in $(z,z')$. We here use the integral reprensentations of the functions $\psi_a$. We conclude the proof in the final step. \\
\paragraph{\textbf{Step 1.}} Assume that $z=N \in \mathbb{N}$ and $z'=N+\beta-1$, and:
\begin{align*}
\tilde{x}:=x+N-1/2, \quad \tilde{y}:=y+N-1/2 \quad \quad \tilde{u}_i \in \mathbb{R}, \quad \tilde{v}_i=\tilde{u}_i+\beta-1, \quad i=1,...,k .
\end{align*}
Then, by means of Propositions \ref{propkernelcharlier} and \ref{propbopsimeixner}, and computations similar to those performed at the beginning of the proof of Theorem \ref{thm1}, formula (\ref{kernelmodifiedzmeas}) makes sense and holds with:
\begin{align*}
\tilde{u}_i=u_i+N-1/2.
\end{align*}
Indeed, the factor $C_{z,z',k}$ comes from the ratio $\frac{c_{N-1}}{c_{N+2k}}$ appearing in Proposition \ref{propkernelcharlier}, while the expressions for $A_k$, $B_k$ and $D_k$ are obtained from Proposition \ref{propkernelcharlier} after performing operations similar to those in the proof of Theorem \ref{thm1}: the polynomials $M_n$ are replaced by the functions $\psi_a$ thanks to homogeneity, and their derivatives $M_n'$ are replaced by $\psi_a'$ by operations on rows. We resume this fact in the following Lemma:
\begin{lem} \label{lemmestep1} Suppose that $z = N \in \mathbb{N}$ is a positive integer, $\beta >0$ is a positive number and that $z' = N +\beta - 1$. Let $\tilde{u}_1,...,\tilde{u}_k \in \mathbb{R}$ be $k$ real numbers, and set:
\begin{align*}
\tilde{v}_i=\tilde{u}_i + \beta - 1, \quad i=1,...,k.
\end{align*}
Let $M_{z,z',\xi}^k$ be the Christoffel deformation of the $z$-measure $M_{z,z',\xi}$ defined by (\ref{defzmeas}), with $u_i$ (resp. $v_i$) replaced by $\tilde{u}_i$ (resp. $\tilde{v}_i$). We then have, for every $n \in \mathbb{N}$ and every $\{\tilde{x}_1,...,\tilde{x}_n \} \subset \mathbb{N}$:
\begin{align*}
M_{z,z',\xi}^k \left( \{x_1,...,x_n \} \subset \mathfrak{S}(\lambda) \right)=M_{z,z',\xi}^k \left( \{\tilde{x}_1,...,\tilde{x}_n \} \subset \mathfrak{S}_N(\lambda) \right)=\det \left( K_{z,z',\xi}^k(x_i,x_j) \right)_{i,j=1}^n
\end{align*}
where $K_{z,z',\xi}^k$ is defined by (\ref{kernelmodifiedzmeas}) and:
\begin{align*}
x_j &= \tilde{x}_j -N +1/2, \quad j=1,...,n
\intertext{and}
u_i &= \tilde{u}_i -N + 1/2, \quad i=1,...,k.
\end{align*}
\end{lem}
\paragraph{\textbf{Step 2. }} For any values of the parameters $z$ and $z'$, and for any $n \in \mathbb{N}$, and any $\{x_1,..,x_n\} \subset \mathbb{Z}'$, Borodin and Olshanski (\cite{bomeixner}, Proposition 3.9) show that the quantity:
\begin{align*}
M_{z,z',\xi} \left( \{x_1,...,x_n \} \subset \mathfrak{S}(\lambda) \right),
\end{align*}
as a function of $\xi$, can be expanded in series around zero, and its coefficients are polynomials in $(z,z')$. As an immediate consequence, we also have that, viewed as a function of $\xi$, the quantity:
\begin{align*}
M_{z,z',\xi}^k \left( \{x_1,...,x_n \} \subset \mathfrak{S}(\lambda) \right)
\end{align*} has an analytic expansion around zero with coefficients that are polynomials in $(z,z')$.\\
\\
\paragraph{\textbf{Step 3}}We will prove that the same statement holds for a slight modification of the kernel $K_{z,z',\xi}(x,y)$. We first assume that $z$ and $z'$ are admissible parameters not in the degenerate series case. Following \cite{bomeixner} (or \cite{bomarkov}), we define:
\begin{align} \label{def:functfzz'}
f_{z,z'}(x)=\frac{\Gamma(x+z'+1/2)}{\sqrt{\Gamma(x+z+1/2)\Gamma(x+z'+1/2)}}
\end{align}
where $x \in \mathbb{Z}'$. By the assumptions on $z$ and $z'$, this function is well defined. Indeed, the Gamma functions are well defined and the fact that:
\begin{align*}
\Gamma(x+z+1/2)\Gamma(x+z'+1/2) >0
\end{align*} allows to define the square root. We also assume that $u_1,...u_k$ are such that $f_{z,z'}(x)$ is well defined for $x \in \{u_1,...,u_k \}$. Observe that interchanging $z$ and $z'$ corresponds to turning $f$ into $1/f$. We will also use the notation:
\begin{align*}
I_a(x;z,z',\xi)=\frac{1}{2i\pi}\oint (1- \sqrt{\xi}\omega)^{-z'+a-1/2} \left(1- \frac{\sqrt{\xi}}{\omega}\right)^{z-a-1/2}\omega^{-x-a}\frac{d\omega}{\omega}
\end{align*}
for any $a \in \mathbb{Z}'$ and $x \in \mathbb{Z}' \cup \{u_1,...u_k \}$. By means of the integral representation (\ref{intreppsi}), we can decompose:
\begin{align} \label{v1}
\psi_a(x;z,z',\xi)=(1-\xi)^{\frac{z-z'+1}{2}}\frac{f_{z,z'}(x)}{f_{z,z'}(-a)}I_a(x;z,z',\xi)
\end{align}
Since the initial definition (\ref{defpsi}) is symmetric in $(z,z')$, we can also write:
\begin{align} \label{v2}
\psi_a(x;z,z',\xi)=(1-\xi)^{\frac{z'-z+1}{2}}\frac{f_{z,z'}(-a)}{f_{z,z'}(x)}I_a(x;z',z,\xi),
\end{align}
Similarly, we have for the derivatives:
\begin{align}\label{v1'}
\psi_a'(x;z,z',\xi)-\frac{f'_{z,z'}(x)}{f_{z,z'}(x)}\psi_a(x;z,z',\xi)=(1-\xi)^{\frac{z-z'+1}{2}} \frac{f_{z,z'}(x)}{f_{z,z'}(-a)}I_a'(x;z,z',\xi)
\end{align}
and:
\begin{align}\label{v2'}
\psi_a'(x;z,z',\xi)+\frac{f'_{z,z'}(x)}{f_{z,z'}(x)}\psi_a(x;z,z',\xi)=(1-\xi)^{\frac{z'-z+1}{2}}\frac{f_{z,z'}(-a)}{f_{z,z'}(x)}I_a'(x;z',z,\xi).
\end{align}
Observe that the integrals $I_a$ and their derivatives $I_a'$ have an analytic expansion in $\sqrt{\xi}$ with coefficients polynomial in $(z,z')$. We thus want to get rid of the factors involving $f_{z,z'}$ and $f_{z',z}$ appearing in formula (\ref{kernelmodifiedzmeas}) through the functions $\psi_a$. We begin with the constant factor $D_k^2$. By (\ref{v1}) and (\ref{v1'}) and their multiplicative structure, we can write 
\begin{multline*}
D_k=(1-\xi)^{k(z'-z+1)}\left(\prod_{i=1}^k f_{z,z'}(u_i)\right)^2 \left( \prod_{i=1}^{2k}f_{z,z'}(-1/2-i+1) \right)^{-1}\\
\times \det\begin{pmatrix}
\left( I_{-1/2-j+1}(u_i;z,z',\xi) \right)_{\substack{i=1,...,k \\ j=1,...,2k}}   \\
\left( I_{-1/2-j+1}'(u_i;z,z',\xi) \right)_{\substack{i=1,...,k \\ j=1,...,2k}}
\end{pmatrix}.
\end{multline*}
But using (\ref{v2}) and (\ref{v2'}), we see that $D_k$ can also be expressed in the following way 
\begin{multline*}
D_k=(1-\xi)^{k(z-z'+1)}\left(\prod_{i=1}^k f_{z,z'}(u_i)\right)^{-2} \left( \prod_{i=1}^{2k}f_{z,z'}(-1/2-i+1) \right)\\
\times \det\begin{pmatrix}
\left( I_{-1/2-j+1}(u_i;z',z,\xi) \right)_{\substack{i=1,...,k \\ j=1,...,2k}}   \\
\left( I_{-1/2-j+1}'(u_i;z',z,\xi) \right)_{\substack{i=1,...,k \\ j=1,...,2k}}
\end{pmatrix}.
\end{multline*}
Thus, writing the determinants $D_k$ in $D_k^2$ in these two different ways leads to:
\begin{align*}
 D_k^2=(1-\xi)^{2k}\det\begin{pmatrix}
\left( I_{-1/2-j+1}(u_i;z,z',\xi) \right)_{\substack{i=1,...,k \\ j=1,...,2k}}   \\
\left( I_{-1/2-j+1}'(u_i;z,z',\xi) \right)_{\substack{i=1,...,k \\ j=1,...,2k}}
\end{pmatrix}
\det\begin{pmatrix}
\left( I_{-1/2-j+1}(u_i;z',z,\xi) \right)_{\substack{i=1,...,k \\ j=1,...,2k}}   \\
\left( I_{-1/2-j+1}'(u_i;z',z,\xi) \right)_{\substack{i=1,...,k \\ j=1,...,2k}}
\end{pmatrix}.
\end{align*}
For the products $B_{k,0}(x)B_{k,1}(y)$ and $B_{k,1}(x)B_{k,0}(y)$, we also write the functions $\psi_a$ and $\psi_a'$ in two different ways, but now, we have to take care of the "boundary factors", i.e. the factors that are not shared by both determinants. For instance, we have 
\begin{multline*}
B_{k,0}(x)B_{k,1}(y)=(1-\xi)^{2k+1}\frac{f_{z,z'}(x)}{f_{z,z'}(y)} \frac{f_{z,z'}(1/2)}{f_{z,z'}(-1/2-2k)} \\
\times \det\begin{pmatrix}
\left( I_{-1/2-j+1}(u_i;z,z',\xi) \right)_{\substack{i=1,...,k \\ j=1,...,2k+1}}   \\
\left( I_{-1/2-j+1}'(u_i;z,z',\xi) \right)_{\substack{i=1,...,k \\ j=1,...,2k+1}}\\
\left(I_{-1/2-j+1}(x;z,z',\xi)\right)_{j=1,..2k+1}
\end{pmatrix}
\det\begin{pmatrix}
\left( I_{1/2-j+1}(u_i;z',z,\xi) \right)_{\substack{i=1,...,k \\ j=1,...,2k+1}}   \\
\left( I_{1/2-j+1}'(u_i;z',z,\xi) \right)_{\substack{i=1,...,k \\ j=1,...,2k+1}}\\
\left(I_{1/2-j+1}(y;z',z,\xi)\right)_{j=1,...,2k+1}
\end{pmatrix},
\end{multline*}
and 
\begin{multline*}
B_{k,0}(y)B_{k,1}(x)=(1-\xi)^{2k+1}\frac{f_{z,z'}(x)}{f_{z,z'}(y)}\frac{f_{z,z'}(-1/2-2k)}{f_{z,z'}(1/2)}\\
\times \det\begin{pmatrix}
\left( I_{1/2-j+1}(u_i;z,z',\xi) \right)_{\substack{i=1,...,k \\ j=1,...,2k+1}}   \\
\left( I_{1/2-j+1}'(u_i;z,z',\xi) \right)_{\substack{i=1,...,k \\ j=1,...,2k+1}}\\
\left(I_{1/2-j+1}(x;z,z',\xi)\right)_{j=1,..2k+1}
\end{pmatrix}
\det\begin{pmatrix}
\left( I_{-1/2-j+1}(u_i;z',z,\xi) \right)_{\substack{i=1,...,k \\ j=1,...,2k+1}}   \\
\left( I_{-1/2-j+1}'(u_i;z',z,\xi) \right)_{\substack{i=1,...,k \\ j=1,...,2k+1}}\\
\left(I_{-1/2-j+1}(y;z',z,\xi)\right)_{j=1,...,2k+1}
\end{pmatrix}.
\end{multline*}
This gives 
\begin{align*}
B&_{k,0}(x)B_{k,1}(y)-B_{k,0}(y)B_{k,1}(x)=(1-\xi)^{2k+1}\frac{f_{z,z'}(x)}{f_{z,z'}(y)}\\
&\times\left( \frac{f_{z,z'}(1/2)}{f_{z,z'}(-1/2-2k)}\det\begin{pmatrix}
\left( I_{-1/2-j+1}(u_i;z,z',\xi) \right)_{\substack{i=1,...,k \\ j=1,...,2k+1}}   \\
\left( I_{-1/2-j+1}'(u_i;z,z',\xi) \right)_{\substack{i=1,...,k \\ j=1,...,2k+1}}\\
\left(I_{-1/2-j+1}(x;z,z',\xi)\right)_{j=1,..2k+1}
\end{pmatrix}
\det\begin{pmatrix}
\left( I_{1/2-j+1}(u_i;z',z,\xi) \right)_{\substack{i=1,...,k \\ j=1,...,2k+1}}   \\
\left( I_{1/2-j+1}'(u_i;z',z,\xi) \right)_{\substack{i=1,...,k \\ j=1,...,2k+1}}\\
\left(I_{1/2-j+1}(y;z',z,\xi)\right)_{j=1,...,2k+1}
\end{pmatrix} \right. \\
&\left. - \frac{f_{z,z'}(-1/2-2k)}{f_{z,z'}(1/2)} \det\begin{pmatrix}
\left( I_{1/2-j+1}(u_i;z,z',\xi) \right)_{\substack{i=1,...,k \\ j=1,...,2k+1}}   \\
\left( I_{1/2-j+1}'(u_i;z,z',\xi) \right)_{\substack{i=1,...,k \\ j=1,...,2k+1}}\\
\left(I_{1/2-j+1}(x;z,z',\xi)\right)_{j=1,..2k+1}
\end{pmatrix}
\det\begin{pmatrix}
\left( I_{-1/2-j+1}(u_i;z',z,\xi) \right)_{\substack{i=1,...,k \\ j=1,...,2k+1}}   \\
\left( I_{-1/2-j+1}'(u_i;z',z,\xi) \right)_{\substack{i=1,...,k \\ j=1,...,2k+1}}\\
\left(I_{-1/2-j+1}(y;z',z,\xi)\right)_{j=1,...,2k+1}
\end{pmatrix} \right) \\
\\
&=\frac{(1-\xi)^{2k+1}}{\sqrt{(z'+1)z'...(z'-2k-1)(z+1)z...(z-2k-1)}} \frac{f_{z,z'}(x)}{f_{z,z'}(y)} \\
&\times \left( (z'+1)z'...(z'-2k-1) \det\begin{pmatrix}
\left( I_{-1/2-j+1}(u_i;z,z',\xi) \right)  \\
\left( I_{-1/2-j+1}'(u_i;z,z',\xi) \right)\\
\left(I_{-1/2-j+1}(x;z,z',\xi)\right)
\end{pmatrix}
\det\begin{pmatrix}
\left( I_{1/2-j+1}(u_i;z',z,\xi) \right)   \\
\left( I_{1/2-j+1}'(u_i;z',z,\xi) \right)\\
\left(I_{1/2-j+1}(y;z',z,\xi)\right)
\end{pmatrix} \right.\\
& \left. - (z+1)z...(z-2k-1) \det\begin{pmatrix}
\left( I_{1/2-j+1}(u_i;z,z',\xi) \right)   \\
\left( I_{1/2-j+1}'(u_i;z,z',\xi) \right)\\
\left(I_{1/2-j+1}(x;z,z',\xi)\right)
\end{pmatrix}
\det\begin{pmatrix}
\left( I_{-1/2-j+1}(u_i;z',z,\xi) \right)   \\
\left( I_{-1/2-j+1}'(u_i;z',z,\xi) \right)\\
\left(I_{-1/2-j+1}(y;z',z,\xi)\right)
\end{pmatrix}  \right),
\end{align*}
where we used:
\begin{align*}
\frac{a}{b}C+\frac{b}{a}D=\frac{a^2C+b^2D}{ab}
\end{align*}
and direct simplifications due to the standard property of the Gamma function 
\begin{align*}
\Gamma(A+1)=A\Gamma(A).
\end{align*} Let us introduce the auxilliary kernel $\tilde{K}_{z,z',\xi}^k(x,y)$:
\begin{align*}
\tilde{K}_{z,z',\xi}^k(x,y)=\frac{f_{z,z'}(y)}{f_{z,z'}(x)}D_k^2 \sqrt{(z'+1)z'...(z'-2k-1)(z+1)z...(z-2k-1)}K_{z,z',\xi}^k(x,y),
\end{align*}
i.e.
\begin{multline*}
\tilde{K}_{z,z',\xi}^k(x,y)\\
= \left( (z'+1)z'...(z'-2k-1) \det\begin{pmatrix}
\left( I_{-1/2-j+1}(u_i;z,z',\xi) \right)  \\
\left( I_{-1/2-j+1}'(u_i;z,z',\xi) \right)\\
\left(I_{-1/2-j+1}(x;z,z',\xi)\right)
\end{pmatrix}
\det\begin{pmatrix}
\left( I_{1/2-j+1}(u_i;z',z,\xi) \right)   \\
\left( I_{1/2-j+1}'(u_i;z',z,\xi) \right)\\
\left(I_{1/2-j+1}(y;z',z,\xi)\right)
\end{pmatrix} \right.\\
\left. - (z+1)z...(z-2k-1) \det\begin{pmatrix}
\left( I_{1/2-j+1}(u_i;z,z',\xi) \right)   \\
\left( I_{1/2-j+1}'(u_i;z,z',\xi) \right)\\
\left(I_{1/2-j+1}(x;z,z',\xi)\right)
\end{pmatrix}
\det\begin{pmatrix}
\left( I_{-1/2-j+1}(u_i;z',z,\xi) \right)   \\
\left( I_{-1/2-j+1}'(u_i;z',z,\xi) \right)\\
\left(I_{-1/2-j+1}(y;z',z,\xi)\right)
\end{pmatrix}  \right),
\end{multline*}
and observe that it has an expansion in $\sqrt{\xi}$ around zero with coefficients that are polynomials in $(z,z')$, since it is the case for the integrals $I_a$ and their derivatives. Observe also that, for any $\{x_1,...,x_n\} \subset \mathbb{Z}'$, we have 
\begin{multline*}
\det\left( \tilde{K}_{z,z',\xi}^k(x_i,x_j) \right)_{i,j=1}^n \\
= D_k^{2n} \sqrt{(z'+1)z'...(z'-2k-1)(z+1)z...(z-2k-1)}^n \det\left( K_{z,z',\xi}^k(x_i,x_j) \right)_{i,j=1}^n.
\end{multline*}
\paragraph{\textbf{Final step. }}Our computations above remain valid for $z$ and $z'$ belonging to the degenerate series case, provided $z > 2k+2$. Assume now that we are in this case, and let $\{x_1,...,x_n\} \subset \mathbb{Z}'$ be a finite subset of $\mathbb{Z}'$. We have established the following equality between series in $\sqrt{\xi}$ with coefficients that are polynomials in $(z,z')$:
\begin{align*}
\left(\det \left( \tilde{K}_{z,z',\xi}^k(x_i,x_j) \right)_{i,j=1,...,n}\right)^2 & = D_k^{4n}\left((z'+1)z'...(z'-2k-1)(z+1)z...(z-2k-1)\right)^n \\
&\times   M_{z,z',\xi}^k \left( \{x_1,...,x_n \} \subset \mathfrak{S}(\lambda) \right)^2 .
\end{align*}
Indeed, the equality follows from the first step and Lemma \ref{lemmestep1}. The right hand side has an analytic expansion in $\xi$ with polynomial coefficients in $(z,z')$, by step 2 and 3, while it also holds for the left hand side by step 3. Since the set:
\begin{align*}\{ (z,z') \in \mathbb{C}^2, z \in \mathbb{N}, \hspace{0.1cm} z > 2k+2 \text{ and } z' \in \mathbb{R}, \hspace{0.1cm} z'>z-1 \}
\end{align*}
is a set of uniqueness for polynomials in two variables, the latter equality holds for all values of $(z,z')$, in particular for all admissible values. This proves the Theorem.
\end{proof}
\subsection{Palm measures of the $z$-measures}
We now state in Theorem \ref{thmpalmzmeas} below the interpretation of Theorem \ref{thmzmeas} in terms of Palm measures of the $z$-measures.
\begin{thm} \label{thmpalmzmeas}
Let $\xi \in (0,1)$ and let $z, z'$ be admissible parameters. Let $u_1, \dots , u_k \in \Z'$ be pairwise distinct and denote by $\P^k_{z,z',\xi}$ the reduced Palm measure of the point process $\mathfrak{S}_*(M_{z,z',\xi})$ at points $u_1,\dots, u_k$. Then, the kernel $K_{z-k,z',\xi}^k$ defined by (\ref{kernelmodifiedzmeas}) is a correlation kernel for $\P^k_{z,z',\xi}$.
\end{thm}
\begin{proof}We start by proving that, up to multiplication by a factor, the correlation functions of the point process $\P^k_{z,z',\xi}$ are analytic in $\sqrt{\xi}$ with coefficients that are polynomials in $(z,z')$. From Proposition \ref{propcvpalm}, the point process $\P^k_{z,z',\xi}$ is a determinantal point process with kernel
\begin{align*}
K^{u_1, \dots , u_k}_{z,z',\xi}(x,y) &= \det \left( K_{z,z',\xi}^0(u_i,u_j) \right)_{1 \leq i,j \leq k }^{-1} \\
&\times \det \begin{pmatrix}
K_{z,z',\xi}^0(x,y) & K_{z,z',\xi}^0(x,u_1) & \dots & K_{z,z',\xi}^0(x,u_k) \\
K_{z,z',\xi}^0(u_1,y) & K_{z,z',\xi}^0(u_1,u_1) & \dots & K_{z,z',\xi}^0(u_1, u_k) \\
\dots & \dots & \dots & \dots \\
K_{z,z',\xi}^0(u_k,y) & \dots & \dots & K_{z,z',\xi}^0(u_k,u_k)     
\end{pmatrix}, 
\end{align*}
where the kernel $K_{z,z',\xi}^0$ is defined by (\ref{kernelmodifiedzmeas}) with $k=0$. By similar arguments as in step 3 in the proof of Theorem \ref{thmzmeas}, the quantity
\begin{align*}
\det \left( K_{z,z',\xi}^0(u_i,u_j) \right)_{1 \leq i,j \leq k }  \cdot K^{u_1, \dots , u_k}_{z,z',\xi}(x,y) 
\end{align*}
is analytic in $\sqrt{\xi}$ with coefficients that are polynomials in $(z,z')$. As an immediate consequence, for any $\{x_1, \dots, x_n \} \subset \Z'$, the quantity
\begin{align*}
\left\{ \det \left( K_{z,z',\xi}^0(u_i,u_j) \right)_{1 \leq i,j \leq k } \right\}^n & \cdot \det \left( K^{u_1, \dots , u_k}_{z,z',\xi}(x_i,x_j) \right)_{i,j= 1}^n \\
&= \left\{ \det \left( K_{z,z',\xi}^0(u_i,u_j) \right)_{1 \leq i,j \leq k } \right\}^n \cdot \P^k_{z,z',\xi} ( \{x_1, \dots, x_n\} \subset X)
\end{align*}
is analytic in $\sqrt{\xi}$ with coefficients that are polynomials in $(z,z')$.

If now $z=N+k$ and $z' \in \R$ with $z' > z-1$, then $\P^k_{z,z',\xi}$ is the Christoffel deformation of the $N$-th Meixner ensemble shifted by $-N +1/2$, i.e.
\begin{align*}
\P^k_{z,z',\xi}( \{x_1,\dots, x_n \} \subset X)= M_{z-k,z',\xi}^k ( \{ x_1,\dots, x_n \} \subset \mathfrak{S}(\lambda)),
\end{align*}
for any $\{x_1,\dots, x_n \} \subset \Z'$. By Lemma \ref{lemmestep1}, we thus have
\begin{align} \label{eq:proofpalmzmeas}
\P^k_{z,z',\xi}( \{x_1,\dots, x_n \} \subset X) = \det \left( K_{z-k,z',\xi}^k(x_i,x_j) \right)_{i,j=1}^n.
\end{align}
Observe now that, up to the multiplication by
\begin{align*}
\left\{\det \left( K_{z,z',\xi}^0(u_i,u_j) \right)_{1 \leq i,j \leq k }\right\}^n,
\end{align*}
both sides of equality (\ref{eq:proofpalmzmeas}) are analytic in $\sqrt{\xi}$ with coefficients that are polynomials in $z,z'$. We deduce that equality (\ref{eq:proofpalmzmeas}) holds for any choice of the parameters $z,z'$, and Theorem \ref{thmpalmzmeas} is proved.
\end{proof}
\subsection{A deformation of the Gamma kernel}
Let $(z,z')$ be admissible parameters, either in the principal or in the complementary series cases. In \cite{bogamma} and \cite{bomarkov}, the authors show that with such a choice of parameters, the $z$-measure, viewed as a point process on $\mathbb{Z}'$, converges as $\xi$ tends to $1$, without any scaling, to a determinantal point process with a correlation kernel expressed only through the Gamma function. They call this kernel the "Gamma kernel". The goal of this section is to show that such a statement still holds for the Christoffel deformation of the $z$-measure of order one, i.e. with $k=1$. Our argument is a slight modification of the one presented in \cite{bogamma}, Theorem 2.3. As mentioned in the introduction, our result must be compared with the one from \cite{bufetovolshanski}, where a hierarchy of Palm deformations of the process with the Gamma kernel is established.
\par
Before stating the Theorem, we give some definitions: for $a \in \mathbb{Z}'$, we define:
\begin{align} \label{def:fungzz'}
g_{z,z'}(a):= \sqrt{\Gamma(z-a+1/2)\Gamma(z'-a+1/2)}\Gamma(-z+a+1/2),
\end{align}
and for $u \in \mathbb{R}$ and $a \in \mathbb{Z}'$, we set:
\begin{align} \label{def:funhzz'}
h_{z,z'}(u,a)=\frac{\Gamma(z'-z)}{f_{z,z'}(u)g_{z,z'}(a)},
\end{align}
recalling that $f_{z,z'}(u)$ has been defined in (\ref{def:functfzz'}). We will denote by $h_{z,z'}'(u,a)$ the derivative of $h_{z,z'}(u,a)$ with respect to $u$. We will also use the notation
\begin{align} \label{def:phi}
\phi_a(u;z,z',\xi)=\psi'_a(u;z,z',\xi)-\frac{\log(\xi)}{2}\psi_a(u;z,z',\xi),
\end{align}
recalling that $\psi_a(u;z;z',\xi)$ has been defined in (\ref{defpsi}).
\begin{thm} \label{thmgamma} Let $(z,z')$ be admissible parameters either in the principal or in the complementary series cases, and let $K_{z,z',\xi}^1$ be the correlation kernel of the Christoffel deformation of order $1$ at $u \in \mathbb{R}$ defined by (\ref{kernelmodifiedzmeas}). We have, for all $x,y \in \mathbb{Z}'$, $x \neq y$:
\begin{align*}
\lim_{\xi \rightarrow 1} (1-\xi)^2 K_{z,z',\xi}^1(x,y)=\frac{\mathcal{C}_{z,z'}}{|(x-u)(y-u)|} \frac{\mathcal{A}_{z,z'}(x,y)}{\Phi_{z,z',0,-1}(u)^2(x-y)}
\end{align*}
where:
\begin{align*}
\mathcal{A}_{z,z'}(x,y)&=\sum_{i,j=0}^2\Phi_{z,z',a(i)b(i)}(u)\Phi_{z,z',a(j)+1,b(j)+1}(u)\Psi_{z,z',i,j}(x,y), \\
\\
\Psi_{z,z',i,j}(x,y)&=h_{z,z'}(x,-1/2+i)h_{z',z}(y,1/2+j)+h_{z',z}(x,-1/2+i)h_{z,z'}(y,1/2+j)\\
&-h_{z,z'}(x,1/2+j)h_{z',z}(y,-1/2+i)-h_{z',z}(x,1/2+j)h_{z,z'}(y,-1/2+i), \\
\\
\Phi_{z,z',a(i),b(i)}(u)&=h_{z,z'}(u,-1/2+a(i))h_{z',z}'(u,-1/2+b(i))\\
&+h_{z',z}(u,-1/2+a(i))h_{z,z'}'(u,1/2+b(i)) \\
&-h_{z,z'}(u,-1/2+b(i))h_{z',z}'(u,-1/2+a(i)) \\
&- h_{z',z}(u,1/2+b(i))h_{z,z'}'(u,-1/2+a(i)),\\
\\
a(0)&=-1, \hspace{0.1cm} a(1) = a(2) =0 \quad; \quad b(0)=b(1)=-2, \hspace{0.1cm} b(2)=-1,\\
\\
\mathcal{C}_{z,z'}&=(z+3)(z+2)(z+1)(z'+2)(z'+1)z'.
\end{align*}
\\
The diagonal entries are obtained by l'Hospital rule.
\end{thm}
\begin{rem}
An explicit expression for the function $\Psi_{z,z',i,j}(x,y)$ may be found from:
\begin{align*}
h_{z,z'}(x,a)h_{z',z}(y,b)&+h_{z',z}(x,a)h_{z,z'}(y,b)-h_{z,z'}(x,b)h_{z',z}(y,a) -h_{z',z}(x,b)h_{z,z'}(y,a)\\
&=\Gamma(z-z')\Gamma(z'-z)\\
&\times\left(\frac{\Gamma(x+z+1/2)\Gamma(x+z'+1/2)\Gamma(y+z+1/2)\Gamma(y+z'+1/2)}{\Gamma(z-a+1/2)\Gamma(z'-a+1/2)\Gamma(z-b+1/2)\Gamma(z'-b+1/2)}\right)^{1/2} \\
&\times \left\lbrace \frac{1}{\Gamma(-z+a+1/2)\Gamma(-z'+b+1/2)\Gamma(x+z'+1/2)\Gamma(y+z+1/2)} \right. \\
&+  \frac{1}{\Gamma(-z'+a+1/2)\Gamma(-z+b+1/2)\Gamma(x+z+1/2)\Gamma(y+z'+1/2)}\\
& - \frac{1}{\Gamma(-z'+a+1/2)\Gamma(-z+b+1/2)\Gamma(x+z'+1/2)\Gamma(y+z+1/2)} \\
&\left. - \frac{1}{\Gamma(-z+a+1/2)\Gamma(-z'+b+1/2)\Gamma(x+z+1/2)\Gamma(y+z'+1/2)}\right\rbrace.\\
\end{align*}
It only involves the Gamma function. A similar expression for the function $\Phi_{z,z',j}$ can be obtained in the same way, and would involve the Gamma function and its first derivative.
\end{rem}
\begin{rem} The factor $(1-\xi)^2$ appearing in our statement should be interpreted as $(1-\xi)^{2k}$ in order to recover Borodin-Olshanski Theorem from \cite{bogamma} and \cite{bomarkov} in the case $k=0$.
\end{rem}
\subsection{Proof of Theorem \ref{thmgamma}}
We begin with the following Proposition which gives the first asymptotic results we need.
\begin{prop} \label{propasypsi1}
For any $u \in \mathbb{R}$, $a \in \mathbb{Z}'$, we have:
\begin{align*}
(1-\xi)^{-\frac{1}{2}}\psi_a(u;z,z',\xi)&=h_{z,z'}(u,a)(1-\xi)^{\frac{z'-z}{2}}\left(1+O(1-\xi)\right)\\
& + h_{z',z}(u,a)(1-\xi)^{\frac{z-z'}{2}} \left( 1 + O(1-\xi) \right) \\
\intertext{and} \\
(1-\xi)^{-\frac{1}{2}}\phi_a(u;z,z',\xi) &= h_{z,z'}'(u,a)(1-\xi)^{\frac{z'-z}{2}}\left(1+O(1-\xi)\right) \\
&+ h_{z',z}'(u,a)(1-\xi)^{\frac{z-z'}{2}} \left( 1 + O(1-\xi) \right)
\end{align*}
as $\xi$ tends to $1$.
\end{prop}
\begin{proof}
This is similar to that in \cite{bogamma}, Theorem 2.3 . We apply the formula:
\begin{align*}
\frac{1}{\Gamma(C)}\left._2F_1\right.(A,B;C;w)&=\frac{\Gamma(B-A)(-w)^{-A}}{\Gamma(B)\Gamma(C-A)}\left._2F_1\right.(A,1-C+A;1-B+A;w^{-1}) \\
&+\frac{\Gamma(A-B)(-w)^{-B}}{\Gamma(A)\Gamma(C-B)}\left._2F_1\right.(B,1-C+B;1-A+B;w^{-1})
\end{align*}
We specialize this results with:
\begin{align*}
A &= -z+a+1/2 \quad ; \quad B=-z'+a+1/2 \quad ; \quad C = u+a+1 \quad ;\\
w&= \frac{\xi}{\xi-1},
\end{align*} 
and plug it into (\ref{defpsi}) to obtain 
\begin{multline*}
(1-\xi)^{-\frac{1}{2}}\psi_a(u;z,z',\xi)\\
=\xi^{\frac{u-a+1}{2}}\left\{h_{z,z'}(u,a)\xi^z(1-\xi)^{\frac{z'-z}{2}}\left._2F_1\right.\left(-z+a+1/2,-z-u+1/2;1+z'-z;\frac{\xi-1}{\xi}\right) \right.\\
\left.+h_{z',z}(u,a)\xi^{z'}(1-\xi)^{\frac{z-z'}{2}}\left._2F_1\right.\left(-z'+a+1/2,-z'-u+1/2;1+z-z';\frac{\xi-1}{\xi}\right) \right\}.
\end{multline*}
By means of the asymptotics:
\begin{align*}
\left._2F_1\right.(A,1-C+A;1-B+A;w^{-1})&=1+O(w^{-1}), \\
\left._2F_1\right.(B,1-C+B;1-A+B;w^{-1})&=1+O(w^{-1}), \\
\frac{d}{dC}\left._2F_1\right.(A,1-C+A;1-B+A;w^{-1})&=O(w^{-1}), \\
\frac{d}{dC}\left._2F_1\right.(B,1-C+B;1-A+B;w^{-1})&=O(w^{-1}),
\end{align*}
which stand for large negative values of $w$ and directly follow from definition (\ref{defgauss}) for the Gauss hypergeometric function, we obtain the desired asymptotics. \end{proof}

These estimates imply the following Lemma.
\begin{lem} \label{lemasypsi} For any $a,b \in \mathbb{Z}+1/2$, and any $x,y \in \mathbb{R}$, we have:
\begin{align*}
\lim_{\xi \rightarrow 1} (1-\xi)^{-1}\left(\psi_{a}(x)\psi_{b}(y)-\psi_{a}(y)\psi_{b}(x) \right)&=h_{z,z'}(x,a)h_{z',z}(y,b)+h_{z',z}(x,a)h_{z,z'}(y,b)\\
&-h_{z,z'}(x,b)h_{z',z}(y,a) -h_{z',z}(x,b)h_{z,z'}(y,a)\\
\intertext{and}
\lim_{\xi \rightarrow 1}(1-\xi)^{-1}(\psi_a(x)\phi_b(x)-\psi_b(x)\phi_a(x))&=h_{z,z'}(x,a)h_{z',z}'(x,b)+h_{z',z}(x,a)h_{z,z'}'(x,b) \\
&-h_{z,z'}(x,b)h_{z',z}'(x,a) - h_{z',z}(x,b)h_{z,z'}'(x,a).
\end{align*}
\end{lem}
\begin{proof}It is similar to that in \cite{bogamma}, Theorem 2.3 and we use the formulas from Proposition \ref{propasypsi1}. Since from (\ref{def:funhzz'}) and (\ref{def:fungzz'}) we have
\begin{align*}h_{z,z'}(x,a)h_{z,z'}(y,b)=h_{z,z'}(x,b)h_{z,z'}(y,a)
\intertext{and}
h_{z',z}(x,a)h_{z',z}(y,b)=h_{z',z}(x,b)h_{z',z}(y,a),
\end{align*}
the terms involving the factors $(1-\xi)^{\pm (z-z')}$ will vanish within a term of the type 
\[(1-\xi)^{\pm (z-z')}O(1-\xi)
\]
in the difference 
\[\psi_a(x)\psi_b(y)-\psi_a(y)\psi_b(x).\]
This term is negligeable because $\mathfrak{R}(z-z') <1$ by our assumption on $z$ and $z'$. Writing down the other terms, we obtain the first part of the Lemma. The proof of the second part is analogous, once we remarked that 
\begin{align*}h_{z,z'}(x,a)h_{z,z'}'(x,b)=h_{z,z'}(x,b)h_{z,z'}'(x,a)
\intertext{and that}
h_{z',z}(x,a)h_{z',z}'(x,b)=h_{z',z}(x,b)h_{z',z}'(x,a).
\end{align*} 
\end{proof}

We can now move to the conclusion of the proof of Theorem \ref{thmgamma}. We consider expression (\ref{kernelmodifiedzmeas}), and first write, as in the proof of Theorem \ref{thm1} 
\begin{multline*}
B_{1,0}(x)B_{1,1}(y)-B_{1,0}(y)B_{1,1}(x)\\
=\sum_{i,j=0}^2 \left[B_{1,0}(x)\right]^i \left[B_{1,1}(y) \right]^j \left( \psi_{-1/2-i}(x)\psi_{1/2-j}(y)-\psi_{-1/2-i}(y)\psi_{1/2-j}(x) \right).
\end{multline*}
By the first part of Lemma \ref{lemasypsi}, we have 
\begin{align} \label{lim1}
\lim_{\xi \rightarrow 1} (1-\xi)^{-1}\left( \psi_{-1/2-i}(x)\psi_{1/2-j}(y)-\psi_{-1/2-i}(y)\psi_{1/2-j}(x) \right)=\Psi_{z,z',i,j}(x,y),
\end{align}
while its second part entails 
\begin{align} \label{lim2}
\lim_{\xi \rightarrow 1} (1-\xi)^{-1}\left[B_{1,p}(x)\right]^i=\lim_{\xi \rightarrow 1} (1-\xi)^{-1}\left[B_{1,p}(y)\right]^i=\Phi_{z,z',a(i)+p,b(i)+p}(u).
\end{align}
Putting (\ref{lim1}) and (\ref{lim2}) together, we obtain 
\begin{align*}
\lim_{\xi \rightarrow 1}(1-\xi)^{-3}\left(B_{1,0}(x)B_{1,1}(y)-B_{1,0}(y)B_{1,1}(x)\right)= \mathcal{A}_{z,z'}(x,y).
\end{align*}
Second part of Lemma \ref{lemasypsi} also implies that 
\begin{align*}
\lim_{\xi \rightarrow 1}(1-\xi)^{-1}D_k = \Phi_{z,z',0,-1}(u),
\end{align*}
which concludes the proof of the Theorem for the off-diagonal part of the kernel. The diagonal entries are treated in a manner similar to that of the proof of Theorem \ref{thm1}, by means of the second part of Lemma \ref{lemasypsi}.


\begin{thebibliography}{99}
\bibitem{andrews} \emph{G.E. Andrews, R. Askey, R. Roy} Special functions {\it Cambridge university press} (2000)
\bibitem{baikdeiftstrahov}\emph{J. Baik, P. Deift, E. Strahov}. Products and ratios of characteristic polynomials of random herimtian matrices {\it Journal of mathematical physics} Vol. 44 No. 8 (2003), 923-975.
\bibitem{bertolaruzza} \emph{M Bertola, G. Ruzza}. Matrix models for stationary Gromov-Witten invariants of the Riemann sphere {\it Preprint} arxiv:2001.10466 (2020) 
\bibitem{borodinbessel} \emph{A. Borodin}. Riemann-Hilbert problem and the discrete Bessel kernel {\it International Mathematics Research Notices} Vol. 2000 (2000)
\bibitem{borodingorin}\emph{A. Borodin, V. Gorin}. Lectures on integrable probability, {\it Probability and Statistical Physics in St. Petersburg, Proceedings of Symposia in Pure Mathematics}, Vol. 91, (2016), 155–214.
\bibitem{boo} \emph{A. Borodin, A. Okounkov, G. Olshanski}. Asymptotics of Plancherel measures for symmetric groups {\it Journal of the american mathematical society} Vol.13 No.3 (2000), 481-515.
\bibitem{bodist}\emph{A. Borodin, G. Olshanski}. Distributions on partitions, point processes, and the hypergeometric kernel {\it Communications in mathematical physics} Vol.211 (2000), 335-358.
\bibitem{bo-zw} \emph{A. Borodin, G. Olshanski}. Harmonic Analysis on the Infinite-Dimensional Unitary Group and Determinantal Point Processes
{ \it Annals of Mathematics, Second Series} Vol. 161, No. 3 ( 2005), 1319-1422.
\bibitem{bogamma}\emph{A. Borodin, G. Olshanski}. Random partitions and the Gamma kernel {\it Advances in mathematics} Vol.194 (2005), 141-202.
\bibitem{bomeixner} \emph{A. Borodin, G. Olshanski}. Meixner polynomials and random partitions {\it Moscow Mathematical Journal} Vol. 6 No.4 (2006), 629-655. 
\bibitem{bomarkov} \emph{A. Borodin, G. Olshanski}. Markov processes on partitions {\it Probability theory and related fields} Vol.135 (2006), 84-152.
\bibitem{giambelli} \emph{A. Borodin, G. Olshanski, E. Strahov} Giambelli compatible point processes {\it Advances in applied mathematics} Vol. 37 No. 2 (2006), 209-248.
\bibitem{borodinstrahov} \emph{A. Borodin and E. Strahov}. Averages of characteristic polynomials in Random Matrix Theory {\it Commun. Pure and Applied Math.} 59 (2006), no. 2, 161-253.
\bibitem{bufetovbessel} \emph{A. I. Bufetov}. A Palm hierarchy for determinantal point processes with the Bessel kernels {\it Proceedings of the Steklov Institute of Mathematics}. Vol. 297 (2017), 90-97.
\bibitem{bufetovquasisymmetries} \emph{A. I. Bufetov} Quasi-symmetries of determinantal point processes {\it Annals of probability}. Vo. 46, No. 2 (2018), 956-1003.
\bibitem{bufetovqiupickrell}\emph{A. I. Bufetov, Y.Qiu}. The explicit formulae for scaling limits in the ergodic decomposition of infinite Pickrell measures {\it Arkiv f\"or matematik}. Vol.54  (2016), 403-435.
\bibitem{bufetovconditional} \emph{A. I. Bufetov}. Conditional measures of determinantal point processes {\it Funct. Anal. Appl.}. Vol. 54 (2020), 7-20.
\bibitem{bufetovlazag} \emph{A. I. Bufetov, P. Lazag} A determinantal point process governed by an integrable projection kernel is Giambelli compatible
 { \it to appear in Annali della Scuola Normale Superiore di Pisa}. (2023). 
\bibitem{bufetovolshanski} \emph{A. I. Bufetov, G. Olshanski}. A hierarchy of Palm measures for determinantal point processes with Gamma kernels {\it Preprint}.  arXiv:1904.13371 (2019)
\bibitem{chhaibinajnudel} \emph{R. Chhaibi, J. Najnudel, A. Nikeghbali}. The circular unitary ensemble and the Riemann
zeta function: the microscopic landscape and a new approach to ratios {\it Inventiones
mathematicæ}, 207 (2017),  23 –113.
\bibitem{daley-verejones} \emph{D.J. Daley, D. Vere-Jones}. An Introduction to the Theory of Point Processes, Vol II, General Theory and Structure {\it Springer Verlag} (2008) 
\bibitem{fyodorovstrahov1} \emph{Y. V. Fyodorov and E. Strahov}. An exact formula for general spectral correlation function of random Hermitian matrices {\it J. Phys.} A 36 (2003) 3203-3213.
\bibitem{ismail}\emph{M.E.H. Ismail}. Classical and quantum orthogonal polynomials in one variable {\it Cambridge university press} (2005)
\bibitem{johanssonplancherel}\emph{K. Johansson}. Discrete orthogonal polynomial ensembles and the Plancherel measure {\it Annals of mathematics} Vol.153 (2001), 259-296.
\bibitem{keatingsnaith} \emph{J. P. Keating, N. C. Snaith}. Random Matrix Theory and zeta(1/2+it) {\it Commun. Math. Phys.} 214 (2000), 57-89.
\bibitem{konig}\emph{W. K\"onig}. Orthogonal polynomial ensembles in probability theory {\it Probability Surveys} Vol. 2 (2005), 385-447.
\bibitem{kontsevich}\emph{M. Kontsevich}. Intersection theory on the moduli space of curves and the matrix Airy function {\it Communications in mathematical physics} Vol. 147 No. 1 (1992), 1-23.
\bibitem{macdonald}\emph{I.G. Macdonald}. Symmetric functions and Hall polynomials {\it Second edition}. {\it Oxford mathematical monographs} (1995)
\bibitem{infinitewedge} \emph{A. Okounkov}. Infinite wedge and random partitions {\it Selecta Mathematica} (2001) 7-57.
\bibitem{ouwong} \emph{Chun-Hua Ou, R. Wong}. The Riemann-Hilbert approach to global asymptotics of discrete orthogonal polynomials with infinite nodes {\it Anal. Appl.}, 8 (2010), 247-286.
http://math.nist.gov/~DLozier/SF21/SF21slides/Ou.pdf
\bibitem{sagan}\emph{B.E. Sagan}. The symmetric group: representations, combinatorial algorithms and symmetric functions {\it Second edition}. {\it Springer} (2001)
\bibitem{shiraitakahashi} \emph{T. Shirai, Y. Takahashi} Random point fields associated with certain Fredholm determinants I: fermion, Poisson and boson point processes. { \it J. Funct. Anal.} 205 (2003), 414-463.
\bibitem{fyodorovstrahov2} \emph{E. Strahov and Y. V. Fyodorov} Universal Results for Correlations of Characteristic Polynomials: Riemann-Hilbert Approach. {\it Commun. Math. Phys.} 241 (2003) 343-382.
\bibitem{szego}\emph{G. Szeg\"o}. Orthogonal polynomials {\it American mathematical society} Vol.23 4th ed. (1973)
\end{thebibliography}
\end{document}